\documentclass[]{interact}
\usepackage{epstopdf}
\usepackage[caption=false]{subfig}
\usepackage[longnamesfirst,sort]{natbib}
\bibpunct[, ]{(}{)}{;}{a}{,}{,}
\renewcommand\bibfont{\fontsize{10}{12}\selectfont}
\usepackage{xcolor}
\usepackage{MnSymbol}
\usepackage{amssymb} \usepackage{amsmath,amscd,mathrsfs}
\usepackage{tikz}
 \usepackage[usenames,dvipsnames]{pstricks}
 \usepackage{epsfig}
 \usepackage{pst-grad} 
 \usepackage{pst-plot} 
 \usepackage[space]{grffile} 

\usepackage{enumitem}
\usepackage{hyperref}
\usepackage[active]{srcltx}

\newcommand{\eop}[1]{
\hspace{10mm} \vspace{-6mm}
\begin{flushright}
$\square_{\text{#1}}$\\ \ \\
\end{flushright}
}
\newenvironment{prueba}[1][{\it Proof}]{\noindent {\it #1.} }{}
\newcommand{\bdem}[1][Proof]{\begin{prueba}[#1]}
\newcommand{\edem}[1][]{\eop{#1}
\end{prueba}}
\def\bsdem{\begin{prueba}[Reference]}\def\bsindem{\begin{proof}[\ ]}

\newtheorem{thm}{Theorem}[section]
\newtheorem{defin}[thm]{Definition}
\newtheorem{remark}[thm]{Remark}
\newtheorem{cor}[thm]{Corollary}
\newtheorem{lem}[thm]{Lemma} \newtheorem{rem}[thm]{Remark}
\newtheorem{prop}[thm]{Proposition} \newtheorem{exam}[thm]{Example}

\newtheorem{notat}[thm]{Notation}
\newtheorem{hecho}[thm]{Fact}
\def\C{\mathbb{C}}
 \def\N{{ \! \rm \ I\!N}}  \def\R{{ \! \rm \ I\!R}}

\def\Z{{ \! \rm Z\!\!Z}}    
      
   \def \square{\hbox {$\sqcup
    $\llap {$\sqcap $}}} 
  
\newcommand{\parcial}[2]{\frac{\partial#1}{\partial#2}}

\begin{document}
\title{A natural semantics for the pullback of fiber bundles of structures}
\author{\name{Leonardo A. Cano G.\textsuperscript{a}\thanks{CONTACT Leonardo A. Cano G. Email: lcanog@unal.edu.co; Pedro H. Zambrano Email: phzambranor@unal.edu.co} and Pedro H. Zambrano\textsuperscript{b}}
\ \\
\affil{\textsuperscript{a,b}Departamento de Matem\'aticas,  Universidad Nacional de Colombia, 
AK 30 $\#$ 45-03 c\'odigo postal 111321, 
Bogot\'a D.C., Colombia.}
}
\maketitle

\begin{abstract}
We remark that forcing on fiber bundles of structures of first order languages is not a compatible semantics with the  pullback (of fiber bundles). Motivated by a combination of epistemology and geometry, we  describe a semantics  which  behaves well with respect to the pullback. This new semantics uses parallel transport in its definition and allows to introduce two different types of extensions for the formulae:  vertical and horizontal extensions.
\end{abstract}


\section{Introduction}
Sheaves of structures on topological spaces correspond to the semantics of {\it Intuitionism} (see~\cite{Caicedo}), located in between of Kripke semantics and topoi logic. This is a paradigm of {\it truth continuity} ({\it continuidad veritativa}, according to~\cite{Caicedo}), which means that if a statement is true in a point therefore it continues being true in a neighborhood of that point. There are other similar approaches of sheaves of structures in several logics (e.g., {\it Continuous Logic} -\cite{OchoaVillaveces}-), where the key idea is still preserving the truth of statements in a neighborhood of a point. This idea was generalized to sheaves based on some special kind of lattices extending the lattice of open-sets of a topology (e.g., locales and quantales, \cite{Elephant,quantales}), which still keeps some geometry behind them and correspond to variants of intuitionism and links topoi and quantum logic.
\\
\\
\indent In this paper, we intend to study this idea  on (smooth) fiber bundles.   Our results appeared as a consequence of playing with a {\it soft} epistemological  interpretation of Ehresmann connections (see Appendix~\ref{sec:background}) which allows to distinguish, for a given proposition, the observer who claims it, the region of space (space--time) where the proposition is claimed, and the accuracy of the measurement on which the proposition holds (see Section~\ref{Sec:epistemological motivation}).  
\\
\\
Examples~\ref{example:No hay pullback forcing}, \ref{exam:incompatibility implication} and \ref{exam: incompatibility one negation} show that the pullback (see Appendix~\ref{sec:background}) is incompatible with respect to classical forcing (see Remark~\ref{Rmk: compatibility Pullback}).  In geometry, the pullback is a very important operation between  fiber bundles. For example, it classifies vector bundles over a given topological space $X$ (let say compact and
Hausdorff), explicitly homotopy classes of continuous functions from $X$ to the Grassmannian
correspond to isomorphic classes of vector bundles (see \cite[Theorem 1.16]{Hatcher}). It also describes
elements of the K--theory of $X$ because they are the pullback of the canonical virtual
class over the Fredholm operators according to Atiyah--Janich theorem (see \cite[Theorem A1]{Atiyah}).
\\
\\
Since the pullback of a fiber bundle is an important geometric operation
and  it is compatible with the notion of fiber bundles of structures (see Definition~\ref{defin: fiber bundle struct} and Proposition~\ref{prop:Pullback of structutures}), it is natural to look for a semantics which is compatible
with  the pullback (in Remark~\ref{Rmk: compatibility Pullback} we are more explicit about what this compatibility means). To do so, we involve our epistemological interpretation of  (Ehresmann) connections into the game (see Section~\ref{Sec:epistemological motivation}), and we introduce the notion of parallel semantics, Definition~\ref{defin:pointwise parallel}. We find out that parallel semantics is compatible with the pullback (see Theorem~\ref{thm:parallel compatible pullback}). In this semantics, the continuity of the truth is defined via curves that play the role of observers moving in space (space--time).
As expected from our epistemological motivation (see Section~\ref{Sec:epistemological motivation}), parallel semantics allows to distinguish three new aspects associated to truth
continuity: space-time stability (truth continuity a la Caicedo), preservation of truth
of statements through the observer movement in space-time (during a time interval)
and stability of the "experimental measure" made by the observer, which we can think as related with the accuracy of the measurements done by the observer. 
\\
\\
To our knowledge, the interaction which we use in this article between Differential Geometry and Mathematical Logic is novel. There have been interaction in other directions. For example, the interaction between Complex Geometry and Model Theory had been explored by multiple authors (see~\cite{MoosaPillay} and references there in), more recently between Differential Geometry and Mathematical Logic (see \cite{Heller-Krol}).
\\ 
\\
In Section~\ref{Sec:epistemological motivation} we present the (epistemological) intuition that lead us to define the parallel semantics in Section~\ref{Sec: parallel semantic}. In Section~\ref{Sec:fibre b structures}, we define fiber bundles of structures and we establish its compatibility with the pullback (see Proposition~\ref{prop:Pullback of structutures}). In Section~\ref{sec:Pointwise semantic}, we define forcing on fiber bundles  in the line of previous   work on sheaves (see~\cite{Caicedo}).  Section~\ref{Sec:compatibility formulae pullback} illustrates through examples and propositions  where in the complexity  of the formula the pullback becomes incompatible with forcing (Remark~\ref{Rmk: compatibility Pullback}), we prove  that the forcing of  a formula without free variables obtained as quantification of an atomic formula is compatible with the pullback (see Proposition~\ref{prop: compatibility qunatified atomic}).  It is in Section~\ref{Sec: parallel semantic}, where the differential structure, through the connection, enters in the game defining what we call parallel semantics.  In Section~\ref{Sec:classic, forcing, pf} using the notion of parallel sections  associated to a connection (see Definition~\ref{Def:parallel section}), we write some results that show how forcing and parallel forcing can be related. 
The fact that the curvature of the connection is $0$ plays an important role in these propositions, so we can  say that the curvature is an obstruction to establish a relation between forcing and parallel forcing (see Corollary~\ref{Cor:Atomic formulas}).  
In Section~\ref{Sec:hor y ver}, we explain how the connection on a fiber bundle allows us to define three different types of extension of a formula: the spatial, the horizontal and the vertical extension. Finally, since this article involves Mathematical Logic and Differential Geometry, which might be considered disconnected branches of mathematics, we include non exhaustive appendixes in both subjects at the end of this article.
\\ \\
\indent Along this article, $M$ and $N$ will denote (smooth) manifolds, $\pi:A\to M$ will denote a fiber bundle and $\mathscr{L}$ will denote a first order language.
 \section{Epistemological motivation}\label{Sec:epistemological motivation}
One of the  epistemological motivations of~\cite{Caicedo} for introducing  forcing on sheaves  is the fact that the subjects  of propositions should be {\it extended} or {\it variables}. This extension or variation of the subjects of propositions is based on the (intuitive) idea that objects and situations  in the world are presented to us as extended in space and time. According to this, for science or conversation, there is no {\it point-wise} subject of propositions, neither {\it instant} phenomenon, because subjects and phenomenons should occupy a 
{\it detectable} region of space and time.  However this notion of extension of the subjects and phenomenons is not defined formally and is left as a soft intuitive motivation in~\cite{Caicedo}.
\\
\\
In this article, we would like to distinguish two {\it extensions} for the subjects of a proposition that were not considered in \cite{Caicedo}. First, we heuristically think as subjects of propositions the measurements of observers and we distinguish for them what we call {\it horizontal and vertical extensions}. We think   the horizontal extension as the region in space (or space--time) where a observer experiments     and the vertical extension as the actual values of these experiments. 
\\
\\
It turns out that using Differential Geometry (explicitly the notion of Ehresmann connection) on  fiber bundles we can distinguish these three epistemological ingredients:
\begin{enumerate}
\item the observer who experiments,
\item the values of her/his experiments and
\item the region of space-time where the experiment is carried out.
\end{enumerate}

Let $\pi:A \to X$ be a fiber bundle with fiber $F$ and let $\Phi$ be an Ehresmann connection (see~Definition~\ref{Def:lift of a path}) on $\pi:A \to X$. 

\begin{center}
\begin{tabular}{|l|l|}
 \hline
 Geometry of fiber bundles&Epistemological interpretation\\
 \hline
 The base space $X$.& Space-time.\\
 The fiber $A_x=\pi^{-1}(x)$.&Experimental measurements.\\
 The connection $\Phi$.&Relates the values of the measurements.\\
 &done at different points of space time.\\
 Paths $\alpha:(-\varepsilon,\varepsilon)\to X$.&Observers moving in space--time.\\
 Horizontal lifts $\tilde{\alpha}:(-\varepsilon,\varepsilon) \to A$.&Measurement $\tilde{\alpha}(t)\in A_{\alpha(t)}$\\
 &at the point $\alpha(t)$ of the space--time.\\
 \hline
 \end{tabular}
 \newline\newline
 Table 1.
\end{center}

In Table~1. we summarize the epistemological interpretation of geometric objects associated to a connection. We interpret $X$ as space--time, the (smooth) paths $\alpha:(-\varepsilon,\varepsilon)\to X$ as observers moving in space--time; given a point $x$ in space--time $X$, we interpret  the fiber $A_x :=\pi^{-1}(x)$ as the possible experimental measurements,  a observer situated at $x$ can make. The connection $\Phi$ is used to relate the values of the measurements done at different points of space time, i.e. if $x_0,x_1 \in X$, then $\Phi$ is used to relate $A_{x_0}$ and $A_{x_1}$. The relation between these fibers is obtained via the  horizontal lifting of paths (see Appendix~\ref{sec:background}). Let $\tilde{\alpha}:(-\varepsilon,\varepsilon) \to A$ be the horizontal lifting  of a path $\alpha:(-\varepsilon,\varepsilon)\to X$ , we interpret that, for the observer $\alpha$, the measurement $\tilde{\alpha}(t)\in A_{\alpha(t)}$ at the point $\alpha(t)$ of the space--time, is equivalent to the measurement  $e_0:=\tilde{\alpha}(0)$ at $x_0:=\alpha(0)$.
\\
\\
This interpretation of the Ehresmann connection has been listened by the first author in discussions about fiber bundles but to our knowledge it is not explicitly written in text books or articles. This article mixes the notion of forcing with this {\it soft}  epistemological interpretation of differential geometry on fiber bundles. As it can be seen, this makes our heuristic   different to the one used by Caicedo in~\cite{Caicedo} even when our intuition can be considered a refinement of Caicedo's point of view.
\section{Fiber bundles of structures and their pullback}\label{Sec:fibre b structures}
\indent Let $\mathscr{L}$ be a first--order signature.  We use the geometric background included in Appendix~\ref{sec:background}.
\begin{defin}\label{defin: fiber bundle struct}{\upshape{(cf. \cite[Definici\'on 2.2]{Caicedo})}}
A fiber bundle $\mathfrak{A}$ of $\mathscr{L}$--structures  is  a  fiber bundle $\pi:A \stackrel{F}{\longrightarrow} M$  such that for each $m \in M$ the fiber $A_m:=\pi^{-1}(m)$ is the universe of an $\mathscr{L}$--structure $\mathfrak{A}_m$  such that
\begin{itemize}
\item[i)] For each relational symbol $R\in \mathscr{L}$ of arity $k<\omega$, the set $R^{\mathfrak{A}} := \bigcup_{m\in M} \{m\} \times R^{\mathfrak{A}_m}$ is an open subset of the direct sum $\bigoplus_{i=1}^k A$ of fiber bundles (see Definition~\ref{DirectSum}).
\item[ii)] For each function symbol $f\in \mathscr{L}$ of arity $k<\omega$, the function $f^\mathfrak{A}:\bigoplus_{i=1}^k A\to  A$ defined by $f^\mathfrak{A}(m,e):=f^{\mathfrak{A}_m}(e) \in A_m$ ($e \in A^k_m$) is a $C^\infty$--function.
\item[iii)]  For each constant symbol $c\in \mathscr{L}$, the function $c^{\mathfrak{A}}:M \to A$ given by $m \mapsto c^{\mathfrak{A}_m}$ is a section of $A$.
\end{itemize}
We will denote this fiber bundle of structures by $\pi:\mathfrak{A}  \stackrel{F}{\longrightarrow} M$.
\end{defin}

The following is a very known fact, which follows from the definition of fiber bundle.

\begin{hecho}\label{existencia-secciones}
For every smooth fiber bundle $\pi:A \to M$ and  every $e\in A$, there exists a local section $s:M \to A$ such that $s(m)=e$.
\end{hecho}

We can define fiber bundles of $\mathscr{L}$--structures of regularity $C^k$ requiring that the sections involved in conditions ii) and iii) given above are $C^k$--sections.
\begin{remark}\label{Rem: identification m cruz Am con Am}
Along the article we will identify $R^{\mathfrak{A}}$ as defined  in item i) of Definition~\ref{defin: fiber bundle struct}), with $ \bigcup_{m\in M}  R^{\mathfrak{A}_m}$.
\end{remark}
The fiber bundle of structures are, via sheafification, examples of sheaves of structures (explained in \cite{Caicedo}).  In a similar way differential manifolds  are examples of   topological spaces but the differential structure allows to define concepts like tangent vectors, tangent bundle or de Rham complex. These concepts simply do not exist on general topological spaces.  In this article we  refine logical or model theoretical concepts on sheaves using smooth fiber bundles, these refinements are not possible on general sheaves of structures without using extra structure.
\begin{exam}
Vector bundles are fiber bundles whose fibers are vector spaces and whose trivializations are linear transformations. All vector bundles over a manifold $M$ are the pullback of the canonical bundle of a Grassmannian   (see  \cite{MardsenTornehabe}). Each vector bundle is an  $\mathscr{L}$--fiber bundle of structures where  $\mathscr{L}$ is the first order signature of $\R$-vector spaces $\{+,\cdot_\alpha :\alpha\in \R \}$. An important example of vector bundle is the tangent space $TM$ of a manifold $M$ (see Remark \ref{Rem:TangetSpace} )
\end{exam}
\begin{exam}
Principal  bundles are very important in Gauge theory, these are fiber bundles whose fibers are groups and whose trivializations are morphism of groups. Each principal bundle is a   $\mathscr{L}$--fiber bundle of structures where  $\mathscr{L}$ is the first order language of groups $\{ \cdot, e, ()^{-1} \}$.
\end{exam}
\begin{exam}
In relativity, space--time is modeled as a $4$--dimensional manifold with a Lorentzian  metric $g$. This metric defines the light cone $C:=g^{-1}((0,\infty))$.   The tangent space $TM$ and the light cone  $C$   conform an example of an $\mathscr{L}$--fibre bundle of structures  for the language of vector spaces and an unary relation symbol $\mathscr{L}:=(+,\cdot_\alpha:\alpha\in \C, R)$. The symbols of sum $+$ and scalar product $\cdot_{\alpha}$ are interpreted as the corresponding sum and scalar products defined in each fiber of $TM$, and $R^{\mathfrak{A}}:=C$.
\end{exam}

\begin{remark}
We observe that given $s_1,\cdots , s_r$ local sections of $A$ and\linebreak
$t(x_1,\cdots x_r)$ an $\mathscr{L}$-term, it is straightforward to see that the fiber bundle function of the fiber bundle $\pi: A \stackrel{F}{\longrightarrow} M$ defined by $m \mapsto t(s_1(m),\cdots s_r(m))$ is in fact a smooth section.
\end{remark}
\begin{defin}
A morphism of fiber bundles  of structures  $\pi_i:\mathfrak{A}_i\stackrel{F_i}{\longrightarrow} M_i$, $i=1,2$ with associated fiber bundles  $\pi_i: A_i\stackrel{F_i}{\longrightarrow} M_i$ is a morphism of fiber bundles which also preserve the $\mathscr{L}$--structure over each fiber. More precisely, a morphism is a pair of smooth maps $\langle \Phi:A_1 \to A_2,\phi: M_1 \to M_2	\rangle$ such that the following diagram commutes$$
\begin{CD}
A_1 @>\Phi>>A_2\\
@VV{\pi_1}V @VV{\pi_2}V\\
M_1 @>\phi>> M_2
\end{CD}
$$
and such that
\begin{itemize}
\item[i)] For each relation symbol $R\in \mathscr{L}$ of arity $k$, if $R^{\mathfrak{A}_1}_m(a_1, \cdots,a_k)$ for $(a_1, \cdots, a_k) \in (A_1)_m$ implies $R^{\mathfrak{A}_2}_{\phi(m)}(\Phi(a_1), \cdots,\Phi(a_k))$.
\item[ii)] For each function symbol $f\in \mathscr{L}$ of arity $k$, $\Phi(f^{(\mathfrak{A}_1)_m}(a_1, \cdots, a_k))=f^{(\mathfrak{A}_2)_{\phi(m)}}(\Phi(a_1), \cdots, \Phi(a_k))$.
\item[iii)] Each constant symbol $c\in \mathscr{L}$ satisfies $\Phi(c^{(\mathfrak{A}_1)_m})=c^{(\mathfrak{A}_2)_{\phi(m)}}$.
\end{itemize}
\end{defin}

The following definition is key in this article,  it shows how to define a fiber bundle of $\mathscr{L}$-structures on a pullback.

\begin{prop}\label{prop:Pullback of structutures}
Let $h:N \to M$ be a smooth function, where $M$ and $N$ are manifolds. Then we can naturally define the pullback $\mathfrak{B}:=h^*(\mathfrak{A})$ of a fiber bundle of $\mathscr{L}$-structures where $\mathfrak{A}$ is a fiber bundle of $\mathscr{L}$-structures  over $M$.
\end{prop}
\bdem
Along this proof, we will use the canonical isomorphisms given in Proposition~\ref{Prop:canonical isomorphisms}. 
Let $\pi:A \to M$ be the fiber bundle that underlies the fiber bundle of structures  $\mathfrak{A}$. We will define a fiber bundle of structures $\mathfrak{B}$ over the fiber bundle $h^*(A)$ as follows:
\begin{itemize}
\item[i)] For each relational symbol $R\in \mathscr{L}$ of arity $k<\omega$, we define
$$R^{\mathfrak{B}} :=\{(n,b_1,\cdots, b_k)\in N\times h^*(A): (\tilde{h}(b_1),\cdots,\tilde{h}(b_k))\in R_{h(n)}^{\mathfrak{A}}\}$$ where $\tilde{h}:h^*(A) \to A$ is the function defined by $\tilde{h}(n,a)=a$, as indicated below Definition~\ref{Def:Pullback fiber bundle}.  Using the identification of  Remark~\ref{Rem: identification m cruz Am con Am}, we notice that $R^{\mathfrak{B}}=\oplus_{i=1}^k h^*(A)\cap \tilde{h}^{-1}(R^\mathfrak{A})$, which proves that $R^{\mathfrak{B}}$ is an open set. 
\item[ii)] For each function symbol $f\in \mathscr{L}$ of arity $k<\omega$, we define the function $f^{\mathfrak{B}}:\bigoplus_{i=1}^k h^*(A)\to  h^*(A)$  by $$f^\mathfrak{B}(n,(n,a_1), \cdots,(n,a_k)):=(n,f^{\mathfrak{A}_{h(n)}}(n,(a_1, \cdots,a_k)) ).$$ By definition of the smooth structure of $\bigoplus_{i=1}^k h^*(A)$, $f^\mathfrak{B}$ is a $C^\infty$--function.
\item[iii)]  For each constant symbol $c\in \mathscr{L}$, we define the function $c^{\mathfrak{B}}:N \to h^*(A)$  by $n \mapsto (n,c^{\mathfrak{A}_{h(n)}})$, which is a section of $h^*(A)$. 
\end{itemize}
\edem[Prop.~\ref{prop:Pullback of structutures}]
\section{Pointwise forcing and local modeling}\label{sec:Pointwise semantic}
Let $\mathfrak{A}$ be a fiber bundle of $\mathscr{L}$--structures. In this section we adapt the forcing of \cite{Caicedo} to fiber bundles.  We recall that in~\cite{Caicedo} sections of the  fiber bundle $\pi:A \stackrel{F}{\longrightarrow} M$ are thought as a kind of nouns of the $\mathscr{L}$--formulae as sentences. 
The main difference between  forcing on sheaves (as explained in~\cite{Caicedo}) and forcing on fiber bundles is that in our new setting we have to impose the locality of the true for equality of terms, because fiber bundles do not have discrete topology in their fibers.

\begin{defin}\label{Def: Pointwise Forcing}{\upshape (cf. \cite[Definici\'on 3.1]{Caicedo})}
Let  $\mathfrak{A}$ be a fiber bundle of $\mathscr{L}$--structures and let $s_1, \cdots, s_r$ be local sections of the fiber bundle  $\pi: A \stackrel{F}{\longrightarrow} M$ defined on a point $m \in M$.   We define recursively on $\mathscr{L}$--formulae, the notion of forcing on the point $m$, for the sections  $s_1, \cdots, s_r$ defined on $m$,
  denoted by  $$
\mathfrak{A} \Vdash_m \varphi (s_1, \cdots, s_r),
$$
as follows:
\begin{itemize}
\item[1)] (atomic case) If $t_1(x_1, \cdots,x_r), \cdots, t_n(x_1, \cdots,x_r)$  are $\mathscr{L}$--terms,
    \begin{itemize} \item[i)]
    	$\mathfrak{A} \Vdash_m (t_1=t_2)(s_1, \cdots, s_r)$ if there exists an open neighborhood  $U\subseteq dom(s_1)\cap \cdots \cap dom(s_k)\subseteq M$ of $m$ such that for all $u \in U$
    $$\mathfrak{A}_{u}\models (t^{\mathfrak{A}_{u}}_1=t_2^{\mathfrak{A}_{u}})(s_1(u), \cdots , s_r(u))$$ for all $u \in U$.
    \item[ii)] If $R \in \mathscr{L}$ is a relational symbol of arity $r<\omega$,  $\mathfrak{A} \Vdash_{m} R(s_1, \cdots, s_r)$ if there exists an open neighborhood  $U\subseteq dom(s_1)\cap \cdots \cap dom(s_k)\subseteq M$ of $m$ such that for all $u \in U$       $$\mathfrak{A}_{u}\models R^{\mathfrak{A}_{u}}(t_1,\cdots t_k)(s_1(u), \cdots , s_r(u)).$$
\end{itemize}
\item[2)] $\mathfrak{A} \Vdash_{m} (\varphi  \land	 \psi)(s_1, s_2, \cdots, s_r)$ if $\mathfrak{A} \Vdash_{m} \varphi (s_1, s_2, \cdots, s_r)$ and $\mathfrak{A} \Vdash_{m}  \psi(s_1, s_2, \cdots, s_r)$.
\item[3)] $\mathfrak{A} \Vdash_{m} (\varphi \vee \psi)(s_1, s_2, \cdots, s_r)$ if $\mathfrak{A} \Vdash_{m} \varphi (s_1, s_2, \cdots, s_r)$ or \linebreak $\mathfrak{A} \Vdash_{m}  \psi(s_1, s_2, \cdots, s_r)$.
\item[4)] $\mathfrak{A} \Vdash_{m} \neg \varphi(s_1, s_2, \cdots, s_r)$ if  there exists an open neighborhood  $U\subseteq M$ of $m$ such that for all $u \in U$ $\mathfrak{A}\not\Vdash_u  \varphi (s_1, \cdots , s_r).$
\item[5)] $\mathfrak{A} \Vdash_{m} (\varphi \to	 \psi)(s_1, s_2, \cdots, s_r)$ if there exists an open neighborhood  $U\subseteq M$ of $m$ such that for all $u \in U$, $\mathfrak{A} \Vdash_{u} \varphi(s_1, s_2, \cdots, s_r)$ implies  $\mathfrak{A} \Vdash_{u} \psi(s_1, s_2, \cdots, s_r)$.
\item[6)] $\mathfrak{A} \Vdash_{m} \exists v \varphi (v,s_1,s_2, \cdots, s_r) $ if there exist a (local) section $s$ defined on $m$ such that $\mathfrak{A} \Vdash_{m}  \varphi (s,s_1,s_2, \cdots, s_r) $.
\item[7)] $\mathfrak{A}  \Vdash_{m} \forall v \varphi (v,s_1,s_2, \cdots, s_r) $ if there exists an open neighborhood  $U\subseteq M$ of $m$ such that for any $u\in U$ and any section $s$  defined on $u$  $\mathfrak{A}  \Vdash_{u}  \varphi (s,s_1,s_2, \cdots, s_r) $
\end{itemize}
\end{defin}
When we say that an atomic formula of relation  is forced in a tuple of sections, we are committing a slight abuse of notation, since formally speaking it should be forced in a section of the direct sum of the fiber bundle (see Definition~\ref{DirectSum}). To solve this abuse of notation we use the canonical isomorphism explained in iii) of Proposition~\ref{Prop:canonical isomorphisms}. The use of this type of identification is usual in geometry without explicitly mentioning the  isomorphism, in this article we try to point out us much as possible the identification used.
\\
\\
In \cite[Definici\'on 3.1]{Caicedo}, it is required  that atomic formulae are true at the point $m$, in contrast to our requirement of being true in an open. It is so because for sheaves, \cite[Lemma 2.2]{Caicedo} guarantees that the lifting of sections for  local homeomorphims implies the stability or extension of atomic formulae involving equalities. 
\begin{defin}\label{Def: local model}
We say that the fiber bundle of $\mathscr{L}$--structures $\mathfrak{A}$ {\bf locally models the $\mathscr{L}$--formula $\varphi$ around $p \in M$ at the sections $s_1,\cdots, s_n$} if there exists an open neighborhood $\emptyset\neq U\subseteq dom(s_1)\cap \cdots \cap dom(s_k)\subseteq M$ of $p$ such that for all $u \in U$, 
$$
\mathfrak{A}_u \models \varphi (s_1 (u),\cdots, s_n(u)). 
$$
\end{defin}
The geometry of the fiber bundle of $\mathscr{L}$--structures makes relational formulae stable or extensive in the sense of the following lemma.
\begin{lem}\label{lemma22 version 1}
{\upshape (cf. \cite[Lemma 2.2]{Caicedo})}
Let $\varphi(x_1,\cdots,x_r)$ be a $\mathscr{L}$--first order formula that contains only the logic operators $\vee, \wedge,\exists$ and atomic formulae without $=$. Let $s_1, \cdots s_r$ be local sections of $A$ defined around a fixed $m\in M$. If $\mathfrak{A}_m\models\varphi(s_1(m), \cdots , s_r(m))$, then there exists some open neighborhood $V$ of $m$ such that $\mathfrak{A}_v\models \varphi(s_1(v), \cdots , s_r(v))$  for all $v \in V$ (i.e. $\mathfrak{A}$ locally models the $\mathscr{L}$--formula $\varphi$ around $m \in M$ at the sections $s_1,\cdots, s_r$).
\end{lem}
\bdem
Let $\overline{x}:=(x_1,\cdots,x_r)$.
\begin{enumerate}
\item Let $R\in \mathscr{L}$ be a relation symbol of arity $k$ and $t_1(\overline{x}),\cdots, t_k(\overline{x})$ be $\mathscr{L}$-terms. Assume that $\mathfrak{A}_m\models R(t_1,\cdots, t_k )(s_1(m),\cdots , s_n(m))$, i.e. $(s_1(m),\cdots , s_n(m))\in (R(t_1,\cdots , t_k))^{\mathfrak{A}_m}$ and therefore\linebreak
$(t_1(s_1(m),\cdots , s_n(m)),\cdots, t_k(s_1(m),\cdots , s_n(m)))\in R^{\mathfrak{A}_m}\subseteq R^{\mathfrak{A}}$. By definition $R^{\mathfrak{A}}$ is an open set of $\oplus_{i=1}^k A$, then there is an open neighborhood $U \subseteq  \oplus_{i=1}^k A$ of $(t_1(s_1(m),\cdots , s_n(m)),\cdots, t_k(s_1(m),\cdots , s_n(m)))$ contained in $R^{\mathfrak{A}}$. Let us denote  $g_i:=t_i^{\mathfrak{A}}(s_1(\cdot),\cdots , s_n(\cdot))$ and $g:=(g_1,\cdots, g_k)$, which is  continuous. Then by continuity of $g$, $V:=g^{-1}(U)$ is an open in $M$ such that for any $m'\in V$ we have that $g(m'):=(t_1(\vec{s}(m')),\cdots ,t_k(\vec{s}(m')))\in U$, where $\vec{s}(\cdot):= (s_1(\cdot),\cdots , s_n(\cdot))$; i.e., $\mathfrak{A}_{m'}\models R(t_1,\cdots, t_k )(\vec{s}(m'))$.
\end{enumerate}
Inductive step:
\begin{itemize}
\item If $\mathfrak{A}_m \models (\varphi_1 \wedge \varphi_2) (s_1(m),\cdots , s_n(m))$  then $\mathfrak{A}_m \models \varphi_1(s_1(m),\cdots , s_n(m))$ and  $\mathfrak{A}_m \models \varphi_2(s_1(m),\cdots , s_n(m))$, and by inductive hypothesis there exist open neighborhoods  $U_1$ and $U_2$ of $m$ such that for all $m' \in U_i$, $\mathfrak{A}_{m'} \models \varphi_i (s_1(m'),\cdots , s_n(m'))$. It is the easy to see that for all $m' \in U_1 \cap U_2$, $\mathfrak{A}_{m'} \models \varphi_1 \wedge \varphi_2 (s_1(m'),\cdots , s_n(m'))$. We can prove the case $\vee$ in an analogous way.
\item Suppose that $\mathfrak{A}_m \models \exists v \varphi(v, s_1(m),\cdots , s_n(m))$, then there is $a \in A_m$ such that $\mathfrak{A}_m \models  \varphi(a, s_1(m),\cdots , s_n(m))$. Notice that there is a local section $s:M \to A$ such that $t(m)=a$.  So we have $\mathfrak{A}_m \models  \varphi(t(m), s_1(m),\cdots , s_n(m))$ and by inductive hypothesis we have  open neighborhood  $U$ of $m$ such that for all $m' \in U$, \linebreak
$\mathfrak{A}_{m'} \models \varphi (s(m'),s_1(m'),\cdots , s_n(m'))$; i.e., $\mathfrak{A}_u \models \exists v\varphi (v,s_1(u),\cdots , s_n(u))$
\end{itemize}
\edem[Lemma~\ref{lemma22 version 1}]
For fiber bundles of $\mathscr{L}$--structures it is easy to provide examples that show that this kind of stability is lost for formulae with equality.
\begin{exam}
Let $\pi_x:\R^2 \to \R$ be a fiber bundle with fiber $\R$ ($\pi_x(x,y)=x$) and consider it as a fiber bundle of $\mathscr{L}$--structures $\mathfrak{A}$ for $\mathscr{L}:=\{=\}$.  Consider the sections $s_1(x)=(x,x)$ and $s_2(x)=(x,-x)$. For  the formula $\varphi(x,y):(x=y)$ we have $\mathfrak{A}_0 \models \varphi(s_1(0),s_2(0))$, but locally $\varphi(s_1,s_2)$ does not hold in $\mathfrak{A}$ around $0$. 
\end{exam}
The following theorem is valid by definition.
\begin{thm}{\upshape (cf. \cite[Teorema 3.1]{Caicedo})}
$\mathfrak{A}  \Vdash_{m}  \varphi (s_1,s_2, \cdots, s_n) $ if and only if there exists and open neighborhood  $U$ of $m$ such that\linebreak
 $\mathfrak{A}  \Vdash_{u}  \varphi (s_1(u),s_2(u), \cdots, s_n(u)) $ for all $u \in U$.
\end{thm}

Notice that pointwise modelling is not equivalent to local modelling, but pointwise forcing is in fact equivalent to local forcing.
\\ \\
\indent The proof of the next theorem follows from Definition ~\ref{Def: Pointwise Forcing} and it is analogous to the proof of~\cite[Teorema 3.2]{Caicedo}. 
\begin{thm}{\upshape (cf. \cite[Teorema 3.2]{Caicedo})}\label{thm:doble negación} Let $\vec{s}:=(s_1,\cdots,s_n)$ be a tuple of local sections defined on $m\in M$. 
$\mathfrak{A}  \Vdash_{m} \neg \neg \varphi (\vec{s}) $ if and only if there exists an open neighborhood  $U$ of $m$ such that  $\{u \in U: \mathfrak{A}  \Vdash_{u} \varphi(\vec{s})\}$ is dense in $U$.
\end{thm}
\section{Compatibility of formulae with the pullback of fiber bundles}\label{Sec:compatibility formulae pullback}
In this section we give examples which illustrate the incompatibility of the pullback (of fiber bundles) and the pointwise forcing. We use the notion of incompatibility with the pullback intuitively but we believe that it would be interesting to formalize it  (to see a synthesis of  the results of this article  around semantics incompatible with pullback see Remark~\ref{Rmk: compatibility Pullback}).  The notions of pullback of a fiber bundle and of a section were reviewed in Appendix~\ref{sec:background}; the pullback of a fiber bundle of structures was explained in Proposition~\ref{prop:Pullback of structutures}.
\\
\\
As explained in Section~\ref{Sec:epistemological motivation}, one can think that a {\it experimental measurement} on a point $m$ in space--time $M$ is given by a tuple $(e_1, \cdots, e_n) \in A^n_m$. For any  $\mathscr{L}$--formula $\varphi$, if $\mathfrak{A}_m \models \varphi(e_1, \cdots, e_n)$, in order to be able of making experimental measures, one would expect that for any observer and what she or he measures, the $\mathscr{L}$--formula $\varphi$ is extensive in time, i.e. $\varphi$ continues being true in some interval of time independently of the movement of the observer. One could wrongly think that this formally means  that if $\mathfrak{A} \Vdash_m \varphi (s_1, \cdots, s_k) $ then for all path (observer) $\sigma:(-1,1) \to M$ such that $\sigma (0)=m$, $\sigma^*(\mathfrak{A}) \Vdash_m \varphi (\tilde{s}_1, \cdots, \tilde{s}_k) $ where $\tilde{s}_i(t):=(t,s_i(\sigma(t)))$ and $\sigma^*(\mathfrak{A})$ is the pullback of $\mathfrak{A}$ as explained in Proposition~\ref{prop:Pullback of structutures}. The following example shows that this is in fact wrong.
\begin{exam}\label{example:No hay pullback forcing}
Let $\mathscr{L}$ be the first order language with a relational symbol $R$ of arity $1$. Let us consider the fiber bundle of $\mathscr{L}$--structures with underlying fiber bundle $\rho:\R^3 \to \R^2$ ($\rho(x,y,z):=(x,y)$) and $R^{\mathfrak{A}}:=\R^3-\{t(1,0,0):t \in \R\}$. Let $s(x,y):=(x,y,0)$ and $\sigma(t)=(t,0)$ ($t\in \mathbb{R}$). Using Theorem~\ref{thm:doble negación}, $\mathfrak{A} \Vdash_{(0,0)} \neg \neg R(s)$ because $\{(x,y) \in \R^2: s(x,y) \in R^{\mathfrak{A}_{(x,y)}}\}$ is dense in $\R^2$. In the same line of reasoning, $\sigma^*\mathfrak{A} \nVdash_{0} \neg \neg R(\sigma^*s)$ where $\sigma^*s(t):=(t,(s\circ\sigma)(t) )=(t,(t,0,0))$, because $\{t \in \R:\sigma^*s(t) \in R^{\sigma^* \mathfrak{A}_{t}} \}=\emptyset$, so it is impossible to find an open neighborhood of $0$ on which this set is dense.
\end{exam}
\begin{remark}\label{rem: incompatibility pullback on sheaves} The previous example can be adapted to sheaves. Let us endow $\R^3=\R^2 \times \R$ with the product topology considering $\R^2$ with the usual topology and $\R$ with the discrete topology. With this topology $\rho:\R^3 \to \R^2$ is a sheaf over $\R^2$ and one can observe that the fiber bundle of structures, the section $s$ and the path $\sigma$ of Example~\ref{example:No hay pullback forcing} provides an example that  forcing is not compatible with the pullback on sheaves. 
\end{remark}
The incompatibility of the pullback points out that {\it the pullback of forcing} is a notion that depends of the observer  (the path or function on which we are taking the pullback). 
\begin{remark}\label{Rmk: compatibility Pullback}
It would be interesting to define the notion of semantics compatible with a pullback in a more general way. For the moment, we have that the classic semantics and local modeling are compatible with the pullback in the sense of Proposition~\ref{Prop:classic semantics compatibiliy pullb} and Corollary~\ref{cor:compatibility local modeling}.    The incompatibility of the pointwise  forcing with the pullback means that there are  formulae which are not compatible with the pullback (as it is proved by Examples~\ref{example:No hay pullback forcing}, \ref{exam:incompatibility implication} and \ref{exam: incompatibility one negation}). 
In  Proposition~\ref{prop: compatibility qunatified atomic} we prove that  the forcing of  a sentence obtained as quantification of an atomic formula is compatible with the pullback.   For parallel semantics (see Section~\ref{Sec: parallel semantic}) the compatibility of the pullback means that Theorem~\ref{thm:parallel compatible pullback} holds.
\end{remark}
\begin{defin}\label{Defin: compatible formulas pullback} Let $\mathscr{L}$ be  a signature.
{\bf A $\mathscr{L}$--formula $\varphi$ is compatible with pullbacks} if for every  fiber bundle of $\mathscr{L}$--structures $\mathfrak{A}$ with underling fiber bundle $\pi:A \to M$, and for every smooth function $f:N \to M$,
\begin{center}
if  $\mathfrak{A} \Vdash_p \varphi (\vec{s})$,  then, for all $q \in f^{-1}[\{p\}]$, $f^* \mathfrak{A} \Vdash_q \varphi (f^*\vec{s})$. 
\end{center}
\end{defin}
Example~\ref{example:No hay pullback forcing} shows that in general no all  formulae are compatible with the pullback. Contrasting, we have the following proposition.

\begin{prop}\label{Prop:classic semantics compatibiliy pullb}
Let $\mathfrak{A}$ be a fiber bundle of structure with underlying fiber bundle $\pi:A \to M$ and let $f:N \to M$ be a smooth function between smooth manifolds. Suppose that $e_1,\cdots,e_k \in A_m$ and $\varphi(x_1,\cdots,x_k)$ is and $\mathscr{L}$--formula. Then,  $\mathfrak{A}_m  \models \varphi (e_1,\cdots e_n)$ if and only if for all $n \in f^{-1}(m)$,  $f^*(\mathfrak{A})_n  \models \varphi ((e_1,n),\cdots (e_k,n))$.
\end{prop}
\bdem
The proof follows by the definition of $\models$ and of pullback of structures (see Proposition~\ref{prop:Pullback of structutures}). 
\edem[Proposition~\ref{Prop:classic semantics compatibiliy pullb}]
As a corollary of the previous proposition we have the following result.

\begin{cor}\label{cor:compatibility local modeling}Let $\mathscr{L}$ be a first order signature.
For local modeling, all the $\mathscr{L}$--formulae are compatible with the pullback; i.e. if $\mathfrak{A}$ is a fiber bundle of structures and $\mathfrak{A}$ locally models $\varphi$ at $p$ for the sections $s_1, \cdots, s_n$, then for all smooth $f:N \to M$ and for all $q \in f^{-1}[\{p\}]$, therefore $f^*\mathfrak{A}$ locally models $\varphi$ at $q$ for the sections $f^* s_1, \cdots, f^* s_n$. 
\end{cor}
\bdem
Let us suppose that $\mathfrak{A}$ locally models $\varphi$ at $p$ for the sections $s_1, \cdots, s_n$. This means that there is a neighborhood $U$ of $p$ such that for all $u \in U$, $\mathfrak{A}_u \models \varphi(s_1(u), \cdots, s_n(u))$.   Let $f:N \to M$ be a smooth function. Since the fibers of  the pullback of  the fiber bundle of structures $\mathfrak{A}$ are {\it essentially} the same fibers of $\mathfrak{A}$ (see Proposition~\ref{prop:Pullback of structutures}), for all $u \in f^{-1}U$, $f^*\mathfrak{A}_u \models \varphi(f^* s_1(u), \cdots, f^*s_n(u))$. \edem[Corollary~\ref{cor:compatibility local modeling}]
\begin{prop}\label{prop: compatibility qunatified atomic}
Let $\mathfrak{A}$ be a fiber bundle of structures for the signature $\mathscr{L}$ with   underlying fiber bundle $\pi:A \to M$. Let $\varphi$ be an $\mathscr{L}$--sentence obtained as quantification of an atomic formula. If $\mathfrak{A} \Vdash_m \varphi$ and $f:N\to M$  is a smooth function, then for all $n\in f^{-1}[\{m\}]$, $f^* \mathfrak{A} \Vdash_n \varphi$.
\end{prop}
\bdem
Let $\varphi$ be an atomic $\mathscr{L}$--formula of the form $\varphi (x_1, \cdots,x_n): h(x_1,\cdots,x_n)=g(x_1,\cdots,x_n)$ where $h$ and $g$ denote terms in the signature $\mathscr{L}$ with free variables $x_1,\cdots,x_n$. Let us suppose $\mathfrak{A} \Vdash_m \forall x_1 \cdots \forall x_n \varphi (x_1, \cdots,x_n)$. Notice that for every  $(a_1,\cdots,a_n) \in (\oplus_{i=1}^n A)_m $ there exists a local section $s$ of $\oplus_{i=1}^n A$  such that  $s(m)=(a_1,\cdots,a_n)$. This fact and  the definition of forcing  $\Vdash_m$ imply that, in fiber bundles, $\mathfrak{A} \Vdash_m \forall x_1 \cdots \forall x_n (h=g) (x_1, \cdots,x_n)$ is equivalent to the fact that $h^{\mathfrak{A}}=g^{\mathfrak{A}}$ in an open subset of $\oplus_{i=1}^n A$ of the form $\pi^{-1}(U)$, where here $\pi$ is the projection $\pi:\oplus_{i=1}^n A \to M$ and $U$ is an open neighborhood of $m$ in $M$. Let us suposse that $U$ is such an open neighborhood  of $m$ for the formula of equality $\varphi$. Then $h^{f^*(\mathfrak{A})}=g^{f^*(\mathfrak{A})}$ in the open subset $\pi^{-1}(f^{-1}(U))$ of $f^*(\oplus_{i=1}^n A)$ where here $\pi$ is the projection $\pi: f^*(\oplus_{i=1}^n A) \to N$.      
\\
\\
Now let us suppose that $\varphi$ is an atomic formula of relation, explictly $\varphi(x_1,\cdots,x_n):R(x_1, \cdots,x_n)$ where $R$ is a symbol of relation in the signature $\mathscr{L}$. By Fact~\ref{existencia-secciones}, there exists a local section $s:M \to A$ such that $s(m)=e$. This implies that  $\mathfrak{A} \Vdash_m \forall x_1 \cdots \forall x_n \varphi (x_1, \cdots,x_n)$ is equivalent to the existence of an open neighborhood $U$ of $m$ in $M$ such that $R^{\mathfrak{A}_u}=\oplus_{i=1}^n A_u $ for all $u \in U$. We observe that if  $R^{\mathfrak{A}_u}=\oplus_{i=1}^n A_u $ for all $u \in U$ then $R^{f^*(\mathfrak{A})_v}=\oplus_{i=1}^n f^*(A)_v $ for $v \in f^{-1}(U)$ because $f^*(A)_v=\{f(v)\} \times A_{f(v )}$ and $R^{f^*(\mathfrak{A})_v}=\{v\} \times R^{\mathfrak{A})_{f(v)}}$. Then, for all $n\in f^{-1}(m)$, $f^* (\mathfrak{A}) \Vdash_n \varphi$. 
\edem[Proposition~\ref{prop: compatibility qunatified atomic}]
\ \\
Next proposition shows that compatible formulae with the pullback are closed under $\land$, $\vee$ and $\exists$.

\begin{prop}\label{prop: compatibility pullback is closed under y,o existe, implica}
Let $\mathfrak{A}$ be a fiber bundle of structures for the signature $\mathscr{L}$ with   underlying fiber bundle $\pi:A \to M$. Let $\varphi$  and $\psi$ be  formulae which are compatible with the pullback. Then, $\varphi \land \psi,$ $\varphi \vee \psi$, $\exists x \varphi$ are also compatible with  pullbacks. 
\end{prop}
\bdem
The proposition is straightforward for $\land$ and $\vee$. Let us suppose that $\mathfrak{A} \Vdash_p \exists x \varphi (\vec{s})$, for some fiber bundle of structures $\mathfrak{A}$ with underling fiber bundle $\pi:A \to M$. Let $f:N \to M$ be a smooth function. By Definition~\ref{Def: Pointwise Forcing}, $\mathfrak{A} \Vdash_p \exists x \varphi (x, \vec{s})$ means that there exists a (local) section $r$ defined in $p$, such that $\mathfrak{A} \Vdash_p  \varphi (r, \vec{s})$. Since $\varphi$ is compatible with the pullback, for all $q \in f^{-1}[\{p\}]$, $f^*\mathfrak{A} \Vdash_q  \varphi (f^*r, f^*\vec{s})$, which is equivalent to $f^*\mathfrak{A} \Vdash_q \exists x \varphi (x, f^*\vec{s})$.  
\edem[Proposition~\ref{prop: compatibility qunatified atomic}]
Next example provides a formula of implication of atomic formulae with equality which is not compatible with the pullback.
\begin{exam}\label{exam:incompatibility implication}
Let $\mathscr{L}$ be a signature with two symbols of functions $f$ and $g$ of arity $1$, and a symbol of constant $0$. Let us consider a fiber bundle of $\mathscr{L}$--structures $\mathfrak{A}$ with underlying fiber bundle $\rho:\R^3 \to \R^2$ ($\rho(x,y,z):=(x,y)$) and such that $$f^{\mathfrak{A}}(x,y,z):=(x,y,9-x^2-y^2), \hspace{1cm} g^{\mathfrak{A}}(x,y,z)=(x,y,4-x^2-y^2)$$ and $0^{\mathfrak{A}_{(x,y)}}:=(x,y,0)$. Let us define $\varphi: f(z)=0$, $\psi: g(z)=0$, and the section  $s(x,y):=(x,y,1)$ of the fiber bundle $\rho:\R^3 \to \R^2$. We have that $$\mathfrak{A} \Vdash_{(3,0)} (\varphi \to \psi)(s)$$ because $\mathfrak{A} \nVdash_{(a,b)} \varphi (s)$  for all $(a,b) \in \R^2$. Let $\sigma:\R \to \R^2$ be the function $\sigma(t):=(3\cos t, 3\sin t)$. We have that $$f^{\sigma^*\mathfrak{A}}(t,3 \cos t,3 \sin t ,z)=(t,0),\hspace{1cm} g^{ \sigma^* \mathfrak{A}}(t,3 \cos t,3 \sin t ,z)=(t,4-9)=(t,-5)$$ and $0^{\sigma^*\mathfrak{A}_{t}}:=(t,0)$.  So, for all $t \in \sigma^{-1}[\{(3,0)\}]=\{2n\pi:n\in\Z\}$,  $\sigma^* \mathfrak{A} \Vdash_{t} \varphi (\sigma^*s)$  but $\sigma^*\mathfrak{A} \nVdash_{t} \psi (\sigma^*s)$, hence   $$\sigma^*\mathfrak{A} \nVdash_t (\varphi \to \psi)(\sigma^* s).$$
\end{exam}

\begin{prop}\label{Prop: negation of relation formulas compatible}
Let $\mathfrak{A}$ be a fiber bundle of structures for the signature $\mathscr{L}$ with   underlying fiber bundle $\pi:A \to M$. Suppose that $R$ is a relational symbol and $\varphi:=\neg R(x_1, \cdots,x_n)$. Then, $\mathfrak{A} \Vdash_p \varphi (\vec{s})$ if and only if there exists an open neighborhood $U$ of $p$ such that:
$$
\{v \in U: \vec{s}(v) \notin R^{\mathfrak{A}_v}\}=U.
$$
\end{prop}
\bdem
From Definition~\ref{Def: Pointwise Forcing} 1) ii) and 4), $\mathfrak{A} \Vdash_p \varphi (\vec{s})$ is equivalent to the existence of  an open neighborhood $U$ of $p$ such that
$$
\{v \in U: \vec{s}(v) \notin R^{\mathfrak{A}_v}\} \text{ is dense in }U.
$$
Since $\{v \in U: \vec{s}(v) \notin R^{\mathfrak{A}_v}\}$ is closed in $U$, because it is the preimage of the open $\bigoplus_{i=1}^n \mathfrak{A} \setminus  R^\mathfrak{A}$ under $\vec{s}$ (which is a continuous function), we are done.
\edem[Proposition~\ref{Prop: negation of relation formulas compatible}]
\ \\
\ \\
The following example shows a negation of an  equality of terms which is incompatible with the pullback.
\begin{exam}\label{exam: incompatibility one negation}
Let $\mathscr{L}$ be a signature with two symbols of function $f$ and $g$ both of arity $1$. Let us consider a fiber bundle of $\mathscr{L}$--structures $\mathfrak{A}$ with underlying fiber bundle $\rho:\R^3 \to \R^2$  defined by $\rho(x,y,z):=(x,y)$ and such that $f^{\mathfrak{A}_{(x,y)}}(z)=(x,y,2x)\text{,  }g^{\mathfrak{A}_{(x,y)}}(z)=(x,y,3x)$. Notice that both $f^{\mathfrak{A}_{(x,y)}}$ and $g^{\mathfrak{A}_{(x,y)}}$ are constant functions. Let us define $\varphi: f(z)=g(z)$ and the section  $s(x,y):=(x,y,1)$ of the fiber bundle $\rho:\R^3 \to \R^2$. We have that 
$$
\mathfrak{A} \Vdash_{(0,0)} \neg \varphi(s)
$$
because for all $(x,y) \in \R^2$,  $\mathfrak{A} \nVdash_{(x,y)}  \varphi(s)$. Consider the path $\sigma(t):=(0,t)$ ($t\in \R$), we have that  
$$\sigma^*\mathfrak{A} \nVdash_{0} \neg  \varphi(\sigma^*s),$$
because for all open neighborhood $U$ of $0$ in $\R$ there is an $u \in U$ such that $\sigma^*\mathfrak{A} \Vdash_{u}   \varphi(\sigma^*s)$.
\end{exam}
The previous examples illustrate  the incompatibility of the forcing with the pullbacks. In  the epistemological interpretation of~\cite{Caicedo} we would say that the sentence $\varphi(x)$ is true in the  extensive subject (the section) $s$ but there would be no interpretation of the incompatibility with the pullback.  In terms of our epistemological interpretation (see Section~\ref{Sec:epistemological motivation},  Table 1), we have a richer interpretation of all the geometric ingredients in the game:
\begin{itemize}
    \item[i)] The section $s$ can be thought as  a way to {\it  horizontally translate  the measurements} done by all observers (i.e. paths of the base space) at a point  in the base space (space-time) at some instant of time. This interpretation of the direct image $s$ is justified since it could be (locally) the integral submanifold  of a Ehresmann connection (thought as a distribution). We recall that we are thinking the parallel translation as the equivalence of measurements in different points of the base space (see Section~\ref{Sec:epistemological motivation}).
    \item[ii)] The path $\sigma$ can be thought as an observer moving in the base space.
    \item[iii)] $\sigma^*(s)$ is the horizontal lift of $\sigma$ under the connection given in the previous item. We can interpret it in the following way: The measurement $\sigma^*(s)(t)$ done at the point $\sigma (t)$ of the base space  would be equivalent to the measurement $\sigma^*(s)(0)$ done at the point $\sigma (0)$ of the base space. Moreover, all measurements $\sigma^*(s)(t)$ are equivalent to each other and we can think of them as the same measurement  performed in different  times and points in the base space.  
    \item[iv)] The sentence $\varphi (x)$ is not true for the  observers (paths) because even when we are evaluating the formula in the same measurement $\sigma^*(s)(t)$, its forcing does not hold, in symbols $\sigma^*(\mathfrak{A}) \nVdash_{0} \varphi (\sigma^*(s)) $.
 \end{itemize}
\section{Parallel forcing of a point} \label{Sec: parallel semantic}
In this section we define parallel forcing, a semantics based in the epistemological motivation given in Section~\ref{Sec:epistemological motivation}. From now on, we will work with a smooth fiber bundle of $\mathscr{L}$--structures  $\mathfrak{A}$ whose fiber bundle
$\pi:A \stackrel{F}{\longrightarrow} M$ has a  connected basis space $M$. We suppose also that $\pi:A \stackrel{F}{\longrightarrow} M$ is endowed with a connection $\Phi$ (see Appendix~\ref{Connection}). 
\begin{defin}\label{defin:pointwise parallel}
Given $e_1, \cdots, e_n \in A_m$ for  a fixed $m \in M$, an $\mathscr{L}$--formula $\varphi(x_1,\cdots,x_n)$ is said to be {\bf $(e_1, \cdots, e_n)$--parallel forced}  if for all path  $\sigma:(-1,1) \to M$ such that $\sigma(0)=m$ we have $\sigma^{*}(\mathfrak{A})\Vdash_{0} \varphi(s_1, \cdots,s_n)$ where the section $s_k$ is the $\sigma^*(\Phi)$--lift to the fiber bundle $\sigma^*(\mathfrak{A})$ of the identity path $i:(-1,1)\to (-1,1)$ such that $s_k(0)=(0,e_k)$ (see Figure 1). We denote it by $\mathfrak{A} ^{\Phi}\Vdash_{e_1, \cdots, e_n} \varphi(x_1, \cdots,x_n)$.
\end{defin}
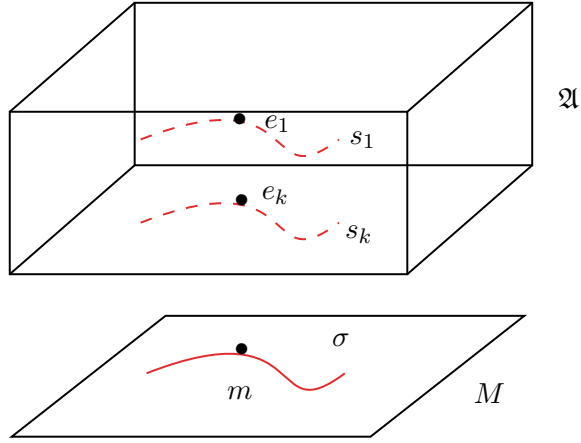
\begin{center}
\tikzset{every picture/.style={line width=0.75pt}} 

\begin{tikzpicture}[x=0.75pt,y=0.75pt,yscale=-1,xscale=1]

\draw    (100, 66) rectangle (300.5, 148)   ;
\draw    (100,66) -- (163,11) ;
\draw    (300.5,66) -- (363.5,11) ;
\draw    (100,148) -- (163,93) ;
\draw    (300.5,148) -- (363.5,93) ;
\draw    (163, 11) rectangle (363.5, 93)   ;
\draw   (178.55,169) -- (359.5,169) -- (281.95,230) -- (101,230) -- cycle ;
\draw [color={rgb, 255:red, 228; green, 42; blue, 42 }  ,draw opacity=1 ]   (169,198) .. controls (253.5,165) and (229,228) .. (269,198) ;

\draw [color={rgb, 255:red, 228; green, 42; blue, 42 }  ,draw opacity=1 ] [dash pattern={on 4.5pt off 4.5pt}]  (166,80) .. controls (250.5,47) and (226,110) .. (266,80) ;
\draw [color={rgb, 255:red, 228; green, 42; blue, 42 }  ,draw opacity=1 ] [dash pattern={on 4.5pt off 4.5pt}]  (166,122) .. controls (250.5,89) and (226,152) .. (266,122) ;
\draw (217,186) node  [align=left] {$\displaystyle \bullet $};
\draw (216,208) node  [align=left] {$\displaystyle m$};
\draw (235,72) node  [align=left] {$\displaystyle e_{1}$};
\draw (234,108) node  [align=left] {$\displaystyle e_{k}$};
\draw (120,141) node  [align=left] {$ $};
\draw (216,70) node  [align=left] {$\displaystyle \bullet $};
\draw (217,111) node  [align=left] {$\displaystyle \bullet $};
\draw (267,181) node  [align=left] {$\displaystyle \sigma $};
\draw (278,80) node  [align=left] {$\displaystyle s_{1}$};
\draw (276,128) node  [align=left] {$\displaystyle s_{k}$};
\draw (342,208) node  [align=left] {$\displaystyle M$};
\draw (382,59) node  [align=left] {$\displaystyle \mathfrak{A}$};
\draw (216,260) node  [align=left] {Figure 1};
\end{tikzpicture}
\end{center}
To explain better the sections $s_k$ in the previous definition, we refer  to Proposition~\ref{Prop: lift of identity as  a section}.
\\
\\
The requirement of two terms being locally equal for the forcing of an equality of terms  (see part 1 i) of~Definition~\ref{Def: Pointwise Forcing} can be considered artificial. Contrasting, the next lemma makes it a consequence of the geometry when we are considering the equality of variables, in analogy to \cite[Lema 2.1]{Caicedo} in the context of sheaves. See Proposition~\ref{prop:atomic equality}  for a generalization of this lemma  for the equality of terms in general.    
\begin{lem}\label{lema caicedo 2.1}
{\upshape (cf. \cite[Lema 2.1]{Caicedo})} If $\tilde{\sigma_1} , \tilde{\sigma}_2: (-\varepsilon,\varepsilon) \to A$ are parallel lifts of a curve $\sigma:(-\varepsilon,\varepsilon) \to M$ associated to some connection $\Phi$ of the fiber bundle $\pi:A \stackrel{F}{\longrightarrow} M$ and for some $t_0 \in (-\varepsilon,\varepsilon)$, $\tilde{\sigma_1}(t_0)= \tilde{\sigma}_2(t_0)$ then there exists a real interval $(t_0-\delta,t_0+\delta)$ such that for all $s \in (t_0-\delta,t_0+\delta)$, $\tilde{\sigma}_1(s)= \tilde{\sigma}_2(s)$.
\end{lem}
\bdem
Recall that the parallel transport of any path $\sigma:(-\varepsilon, \varepsilon) \to M$ is a solution of the equation
$$
\frac{d}{dt}\tilde{\sigma}(t)=d \pi^{-1}(\sigma'(t)),
$$
where we are considering $d \pi$ as an isomorphism  of the horizontal space $HA$ (induced by $\Phi$) and $TM$. This Lemma is a direct consequence of the theorem of existence and uniqueness of solutions of ordinary differential equations. \edem[Lemma~\ref{lema caicedo 2.1}]

Parallel forcing is preserved under isomorphisms of fiber bundles with connections (see Definition~\ref{defin: isomorphic bundles with connection}). 
It implies the following lemma  which is  our main tool to prove that parallel semantics (see Definition~\ref{defin:pointwise parallel}) is in fact compatible with the pullback of any smooth function, not only with respect to paths. 

\begin{lem}\label{lema:parallel compatible pullback}
Let $f:N \to M$ and $g:M \to P$ be  smooth functions ($M$, $N$ and $P$ manifolds), let $\pi:A \stackrel{F}{\to} P$ be a fiber bundle. Suppose  that $\mathfrak{A}$ is a fiber bundle of structures with underling fiber bundle $\pi:A \stackrel{F}{\to} P$ over the signature $\mathscr{L}$. Then, for all formula $\varphi$ in the signature $\mathscr{L}$, and all sections $s_1,\cdots, s_k$ of $\pi:A \to P$, and all $n_0 \in N$, we have:
\begin{center}
$f^* (g^* \mathfrak{A} )   \Vdash_{n_0} \varphi (s_1,\cdots, s_k)$ if and only if  $(g \circ f)^* \mathfrak{A} \Vdash_{n_0} \varphi (\tilde{s}_1,\cdots, \tilde{s}_k)$,
\end{center}
where $\tilde{s}_i$ is the section of ($g \circ f)^* A$ induced by the section $s_i(n)=(n,f(n),T_i(n))$ of $f^*(g^*A)$ through the canonical isomorphism between $f^*(g^*A)$  and  $(g \circ f)^* A$ defined in Proposition~\ref{Prop:canonical isomorphisms}, specifically $\tilde{s}_i(n)=(n,T_i(n))$.
\end{lem}
\bdem
Doing induction on formulas one can observe that the forcing at $n \in N$ of the formula $\varphi$, on the fiber bundles of structures $(g \circ f)^* \mathfrak{A}$ and $f^*(g^*\mathfrak{A})$ for the sections $s_1,\cdots, s_k$ and $\tilde{s}_1,\cdots, \tilde{s}_k$ respectively, depends in the same way of the functions $T_i: N \to A$. 
\edem[Lemma~\ref{lema:parallel compatible pullback}]

\begin{thm}\label{thm:parallel compatible pullback}
Let $f:N \to M$ be a smooth function and let $\pi:\mathfrak{A} \to M$ be a fiber bundle of structures. If $e_1, \cdots,e_k \in A_m$ then,
\begin{center}
If $\mathfrak{A}^\Phi \Vdash_{e_1, \cdots,e_k} \varphi$, then for all $n \in f^{-1}(m)$, $f^*(\mathfrak{A})^{f^*(\Phi)} \Vdash_{(n,e_1), \cdots,(n,e_k)} \varphi$.
\end{center}
\end{thm}
\bdem
Let us suppose $\mathfrak{A}^\Phi \Vdash_{e_1, \cdots,e_k} \varphi$. Let $\sigma:(-\epsilon, \epsilon) \to N$ be such that $\sigma(0)=n$ for $n \in f^{-1}[ \{m\}]$. According to Lemma~\ref{lema:parallel compatible pullback},   
\begin{center}
$(f \circ \sigma)^* \mathfrak{A} \Vdash_{0} \varphi (s_1,\cdots, s_k)$ if and only if $\sigma^* (f^* \mathfrak{A} )   \Vdash_{0} \varphi (\tilde{s}_1,\cdots, \tilde{s}_k)$,
\end{center}
for $s_i$ the section induced by the $(f \circ \sigma)^*(\Phi)$--horizontal lift of the identity of $\R$ to the fiber bundle $(f \circ \sigma)^* A$ such that $s_i(0)=(0,e_i)$ and $\tilde{s}_i$ is the section of  $\sigma^* (f^*A)$ defined by  $\tilde{s}_i(t)=(t,\sigma(t),r_i(t))$ where $r_i(t)$ is the function such that   $s_i(t)=(t,r_i(t))$. Using Propositions~\ref{prop: horizontal lift pullback} and~\ref{Prop:canonical isomorphisms}, we can see that  $\tilde{s}_i$ is equal to the $\sigma^* (f^* \Phi)$--horizontal lift of the identity of $\R$ to the fiber bundle $\sigma^* (f^* A )$.
\edem[Theorem~\ref{thm:parallel compatible pullback}]
\section{Classic semantics, forcing and parallel forcing} \label{Sec:classic, forcing, pf}
In this section we illustrate differences and similarities between forcing, parallel forcing and the classical semantics. We continue also clarifying  the epistemological motivation given in Section~\ref{Sec:epistemological motivation}.     
\\
\\
Examples~\ref{example:No hay pullback forcing},\ref{exam:incompatibility implication} and \ref{exam: incompatibility one negation} distinguish forcing (Definition~\ref{Def: Pointwise Forcing}) from  parallel forcing by proving that the former one is not compatible with pullbacks.  We proved that  parallel forcing is compatible with pullbacks in Theorem~\ref{thm:parallel compatible pullback}. The classic semantics ($\models$) is also compatible with pullbacks (see Proposition~\ref{Prop:classic semantics compatibiliy pullb}) which also distinguishes it from parallel forcing.  Next example shows a difference between  classic semantics and parallel forcing. 
\begin{exam}\label{Example:Classic versus parallel}
On the fiber bundle $\rho:\R^3 \to \R^2$ ($\rho(x,y,z)=(x,y)$) we consider the connection $\Phi=dz \otimes \parcial{}{z}$ (see Example~\ref{Exam:trivial fiber bundle with connect}).  Let $\mathscr{L} $ be a  signature with a  relational symbol $R$  of arity $1$. Let us suppose that $\mathfrak{A}$ is a fiber bundle of structures with underling fiber bundle $\rho:\R^3 \to \R^2$ and that $R^{\mathfrak{A}}:=\{(x,y,z)\in \R^3: x-2y+z \neq 0 \}$, in other words in each fiber we have $R^{\mathfrak{A}_{(a,b)}}=\{(x,y,z):x=a, y=b, a-2b+z \neq 0\}$ . Then, we have $\mathfrak{A}_{(0,0)} \models \neg R(0,0,0)$ but at $e_1:=(0,0,0)\in A_{(0,0)}$ we have that $\mathfrak{A} \nVdash^\Phi_{e_1} \neg R(x) $ because for $\sigma(t)=(t,t)$, the $\Phi$-horizontal lift  $\tilde{\sigma}$ such that $\tilde{\sigma}(0)=(0,0,0)$ is given by  $\tilde{\sigma}(t)=(t,t,0)$, but for every neighborhood of $t=0$, there exists an $s$ such that $\tilde{\sigma}(s)$ belong to $R^{\mathfrak{A}}$.
\end{exam}
Next we study  relations between the classic semantics, forcing and the parallel forcing with respect to atomic formulas of equality. We have the following  characterization of  parallel forcing for equality of terms.
\begin{prop}\label{Prop:Phi forcing equality more natural}
Let $\mathfrak{A}$ be a fiber bundle of structures and let $\Phi$ be a connection on $\pi:A \to M$, the underlying fiber bundle. Suppose that $f$ and $g$ are $\mathscr{L}$--terms. Then  the following are equivalent:
\begin{itemize}
    \item[i)] $\mathfrak{A} \Vdash^\Phi_{(e_1,\cdots,  e_l)} (f=g)(x_1,\cdots , x_l)$, where $e_1,\cdots , e_l \in A_m$.
    \item[ii)] For all $\sigma:(-\varepsilon, \varepsilon) \to M$ such that $\sigma(0)=m$, there exists a $\delta>0$ such that  $\mathfrak{A}_{\sigma(t)} \models (f=g)(\tilde{\sigma}_1(t),\cdots \tilde{\sigma}_l(t))$ for $t \in (-\delta,\delta )$ where $\tilde{\sigma}_i$ is the $\Phi$--lift of $\sigma$ such that $\tilde{\sigma}_i(0)=e_i$.
    \item[iii)] For all $\sigma:(-\varepsilon, \varepsilon) \to M$ such that $\sigma(0)=m$, there exists a $\delta>0$ such that $\mathfrak{A}_{\sigma(t)} \models (f=g)(\tilde{\sigma}(t))$ for $t \in (-\delta,\delta )$, where $\tilde{\sigma}$ is the $\oplus_{i=1}^l \Phi$--lift of $\sigma$ such that $\tilde{\sigma}(0)=(e_1,\cdots e_l)$.
\end{itemize}
\end{prop}
\bdem
The equivalence of i) and ii) follows from the observation that if $s_i$ denotes the $\sigma^*(\Phi)$--lift to $\sigma^*(A)$ of the identity $i:\R\to \R$, such that $s_i(0)=(0,e_i)$, we have  $s_i(t)=(t, \tilde{\sigma}_i(t))$.  
\\
\\
ii) and iii) are equivalent due to Proposition~\ref{Prop:lift of Ak}.
\edem[Proposition~\ref{Prop:Phi forcing equality more natural}]

The following Corollary says that the parallel forcing implies the classical pointwise semantics.

\begin{cor}\label{Cor:forcing equality implies classic}
Let $\mathfrak{A}$ be a fiber bundle of structures and let $\Phi$ be a connection on $\pi:A \to M$ the underlying fiber bundle. Let $f$ and $g$ be terms of the signature $\mathscr{L}$, and let $e_1,\cdots,e_n\in A_m$ for $m \in M$. If $\mathfrak{A} \Vdash_{(e_1, \cdots,e_n)}^\Phi (f=g)(x_1,\cdots,x_n)$,\linebreak then $\mathfrak{A}_m \models (f =g)(e_1, \cdots,e_n)$.
\end{cor}
\bdem
This follows directly from the definition of $\Vdash^\Phi_{(e_1, \cdots,e_n)}$.
\edem[Corollary~\ref{Cor:forcing equality implies classic}]

Notice that we also have the converse of the previous Corollary for variables as $\mathscr{L}$--terms.

\begin{prop}\label{Prop:eq variables forcing y p.forcing}
Let $\mathfrak{A}$ be a fiber bundle of structures and let $\Phi$ be a connection on $\pi:A \to M$ the underlying fiber bundle. Suppose that $e_1,e_2 \in A_m$ for $m\in A$, then, $\mathfrak{A}_m \models (x=y)(e_1,e_2)$ if and only if $\mathfrak{A}\Vdash_{e_1,e_2}^\Phi (x=y)$.
\end{prop}
\bdem
It is a direct consequence of Lemma~\ref{lema caicedo 2.1} and Corollary~\ref{Cor:forcing equality implies classic}. \edem[Proposition~\ref{Prop:eq variables forcing y p.forcing}]  

In contrast with Proposition~\ref{Prop:eq variables forcing y p.forcing}, next example shows that we do not have the previous equivalence for forcing (see Definition~\ref{Def: Pointwise Forcing}).


\begin{exam}
Consider the fiber bundle $\rho:\R^3 \to \R^2$ ($\rho(a,b,c)=(a,b)$). We have $\mathfrak{A}_{(0,0)} \models (x=y)((0,0,1),(0,0,1) ) $, but we can construct sections $s_1$ and $s_2$  for which $s_1(0,0)=(0,0,1)=s_2(0,0)$ and   $\mathfrak{A}_{(0,0)} \nVdash_m (x=y)(s_1,s_2) $. For instance take $s_1(x,y)=(x,y,1+x)$ and $s_2(x,y)=(x,y,1+y)$. 
\end{exam}
Next example shows that the converse of Corollary~\ref{Cor:forcing equality implies classic} is not true for equality of terms, i.e. there are equalities of terms which are true in a fiber but which are not parallel forced. 
\begin{exam}\label{exam: equality of terms in fiber does not imply forcing of it}
Let $\rho:\R^3 \to \R^2$ be the fiber bundle with connection $\Phi$ of Example~\ref{Example:Classic versus parallel}. Let $\mathscr{L} $ be a  signature with  symbols of function $f$ and $g$ of arity one. Let $\mathfrak{A}$ be a fiber bundle of structures on which $f^\mathfrak{A}(a,b,c)=(a,b,c+a)$ and $g^\mathfrak{A}(a,b,c)=(a,b,c+b)$. 
Clearly,   $\mathfrak{A}_{(0,0)}\models  (f(x)=g(x))(0,0,0)$. Let us consider the path $\sigma(t):=(t,2t)$. 
Example~\ref{Exam:trivial fiber bundle with connect}, shows that $\tilde{\sigma}(t)=(t,2t,0)$, where $\tilde{\sigma}$ is the $\Phi$--horizontal lift of the path $\sigma$ such that $\tilde{\sigma}(0)=(0,0,0)$. 
 If we denote by $p$ and $q$ the interpretation of the symbols of function  $f$ and $g$ in $\sigma^* \mathfrak{A}$ (see Proposition~\ref{prop:Pullback of structutures}), we have
$f^{\sigma^*\mathfrak{A}}(t,(t,2t,v)):=(t,t,2t,v+t)$ and $g^{\sigma^*\mathfrak{A}}(t,(t,2t,v)):=(t,t,2t,v+2t)$, then $\sigma^* \mathfrak{A} \nVdash_0  (f(x)=g(x))(s) $ where $s$ is $\sigma^*\Phi$--horizontal lift of the identity. Hence $\mathfrak{A}^\Phi \nVdash_{(0,0,0)} (f(x)=g(x))$.
\end{exam}
In terms of the epistemological motivation of Section~\ref{Sec:epistemological motivation}, the previous example shows that  the observation of equality of two terms at a point is not sufficiently stable, it depends on where the observer is located and it might depend of a very particular observer, the one who stays at the point where the observation is made. The formula $(f(x)=g(x))$ in Example~\ref{exam: equality of terms in fiber does not imply forcing of it} is not parallel forced because of this lack of stability and this dependence of the observer.
\\
\\
 Next proposition proves that for fiber bundles with connections $\Phi$ of curvature $0$, $\Phi$--parallel forcing is  equivalent to force on {\it naturally} associated sections. The main  tool to prove this is that for connections of curvature $0$ (see Definition~\ref{curvature}) the horizontal bundle is an integrable distribution and hence it induces a local section (see Proposition~\ref{prop: curvature 0 sections}).
\begin{prop}\label{prop:atomic equality}Let $\mathfrak{A}$ be a fiber bundle of structures and let $\Phi$ be a connection on $\pi:A \stackrel{F}{\longrightarrow} M$ of curvature $0$ on the underlying fiber bundle. Suppose that $f$ and $g$ are $\mathscr{L}$--terms. Then 
$\mathfrak{A} \Vdash^\Phi_{(e_1,\cdots, e_l)} (f=g)(x_1,\cdots, x_l)$ if and only if $\mathfrak{A} \Vdash_{m} (f=g)(s_1,\cdots, s_l)$
where  $s_i$ is a    $\Phi$--parallel (local) sections of   $A$  such that $s_i(m)=e_i$ defined on a sufficiently small open set of $M$. 
\end{prop}
\bdem
Since $\Phi$ has curvature $0$, we have for each $e_i \in A_m$, there exists $U_i \subset M$ an open neighborhood of $m$ and $s_i:U_i \to A$ such that $s_i$ is a $\Phi$--lift of the inclusion $j:U_i \to M$ and $s_i(m)=e_i$ (see Proposition~\ref{prop: curvature 0 sections}).  This is equivalent to note that each $s_i$ is the integral subvariety of the horizontal bundle $HA_\Phi$ at $e_i$.   Equation (\ref{eq: local section and parallel lift}) and Proposition~\ref{Prop:Phi forcing equality more natural} prove that if $\mathfrak{A} \Vdash_{m} (f=g)(s_1,\cdots s_l)$, then $\mathfrak{A} \Vdash^\Phi_{(e_1,\cdots, e_l)} (f=g)(x_1,\cdots x_l)$.  
\\
\\
We use  ii) of  Proposition~\ref{Prop:Phi forcing equality more natural} to prove   that there exists an open neighborhood $U\subset M$ of $m$ such that  for all $n\in U$, 
$$\mathfrak{A}_{n} \models (f=g)(s_1(n),\cdots, s_k(n)).$$
We proceed by contradiction. Suppose that we can not find such neighborhood $U$, then there exists a sequence $(a_l)_{l \in \N}\subset M$  which converges to $m$ (as fast as needed) and such that $$\mathfrak{A}_{a_l} \nmodels (f=g)(s_1(a_l),\cdots s_k(a_l)).$$ By interpolation methods,  there exists a path $\sigma:(-\varepsilon,\varepsilon) \to M$ such that  $\{a_k: k \in \N\}$ is contained in  $\sigma(-\varepsilon,\varepsilon)$. Since by Equation (\ref{eq: local section and parallel lift}) the $\Phi$--lift at $m$ of $\sigma$ is given by $\tilde{\sigma}=s\circ \sigma$, we have that Proposition~\ref{Prop:Phi forcing equality more natural} ii) does not hold, which is a contradiction. This finishes  the proof.  \edem[Proposition~\ref{prop:atomic equality}]
\begin{cor}\label{Cor:Atomic formulas}
Let $\mathfrak{A}$ be a fiber bundle of structures and let $\Phi$ be a connection on $\pi:A \to M$, the underlying fiber bundle. Suppose that the curvature of $\Phi $ is $0$. Then for every atomic formula $\varphi(x_1,\cdots,x_k)$ we have that the following are equivalent:
\begin{itemize}
\item[i)]$\mathfrak{A} \Vdash^\Phi_{(e_1,\cdots, e_k)} \varphi(x_1,\cdots,x_k)$.
\item[ii)] There exist $s_1, \cdots,s_k$ $\Phi$--parallel (local) sections of   $A$ whose domains contain $m$   such that   $\mathfrak{A} \Vdash_{m} \varphi(s_1,\cdots,s_k)$. 
\end{itemize}
\end{cor}
\bdem
For equality of $\mathscr{L}$-terms, it follows from Proposition~\ref{prop:atomic equality}. 
If $\varphi=R(f_1,\cdots, f_l)(x_1,\cdots,x_k)$ where each $f_i$ is a $\mathscr{L}$-term and $\mathfrak{A} \Vdash^\Phi_{(e_1,\cdots, e_k)} \varphi(x_1,\cdots,x_k)$, and therefore $\mathfrak{A}_m \models \varphi(e_1,\cdots, e_k)$,
then $(f^{\mathfrak{A}_m}_1,\cdots, f^{\mathfrak{A}_m}_l)(e_1,\cdots e_k)$ belongs to $R^{\mathfrak{A}_m}$ and  there exists an open neighborhood $\tilde{U}$ of $(f^{\mathfrak{A}_m}_1,\cdots, f^{\mathfrak{A}_m}_l)(e_1,\cdots e_k)$ in $\oplus_{i=1}^l A$ such that $\tilde{U} \subset R^{\mathfrak{A}}$. Let us denote $U:=(f^{\mathfrak{A}}_1,\cdots, f^{\mathfrak{A}}_l)^{-1}(\tilde{U})$. $U$ is open in $\oplus_{i=1}^k A$, the domain of $(f^{\mathfrak{A}}_1,\cdots, f^{\mathfrak{A}}_l)$. Hence for arbitrary local sections $s_1, \cdots, s_k$ of $A$ such that for all $m \in \pi^{-1}(U)\subset M$, $(s_1, \cdots, s_k)(m) \in U$, we have $\mathfrak{A} \Vdash_{m} \varphi(s_1,\cdots,s_k)$. From this fact follows the proposition.
\\
\\
Let us suppose that $\mathfrak{A} \Vdash_{m} \varphi(s_1,\cdots,s_k)$ where  $s_1, \cdots,s_k$ are $\Phi$--parallel (local) sections. Since $\varphi$ is atomic, $\mathfrak{A} \Vdash_{m} \varphi(s_1,\cdots,s_k)$ is equivalent to the existence of an open neighborhood $U$ of $m$ in $M$ such that $\mathfrak{A}_u \models \varphi(s_1(u),\cdots,s_k(u))$ for all $u \in U$. This and Equation (\ref{eq: local section and parallel lift}) imply that  $\mathfrak{A}_{\sigma(t)} \models \varphi(\tilde{\sigma}_1(t),\cdots,\tilde{\sigma}_k(t))$, and  the corollary follows from the definition of parallel forcing. \edem[Corollary~\ref{Cor:Atomic formulas}] 
\begin{rem}
Notice that the  hypothesis of curvature $0$ is  needed in the previous Corollary~\ref{Cor:Atomic formulas} for existence of  $\Phi$--parallel (local) sections (see Definition~\ref{Def:parallel section} and Proposition~\ref{prop: curvature 0 sections}).
\end{rem}
\begin{rem}
Corollary~\ref{Cor:Atomic formulas} can be extended  to  formulae constructed from atomic formulas using disjunction and conjunction.
\end{rem}
 However the next example shows that  forcing of a formula on parallel sections (see Definition~\ref{Def:parallel section}) is not equivalent to parallel forcing. 
\begin{exam}
Let $\mathscr{L}$ be a signature with a unary symbol of relation $R$ and let $\mathfrak{A}$ be a fiber bundle of structures (see Definition~\ref{defin: fiber bundle struct}) with underling fiber bundle $\pi: \R^3 \stackrel{\R}{\to} \R^2$ ($\pi(x,y,z)=(x,y)$) and such that $R^{\mathfrak{A}}:=\R^3-\{(x,y,z) \in \R^3: x=y, z=0\}$. Let $\Phi:=\parcial{}{z} \otimes dz$ be a connection on $\pi: \R^3 \stackrel{\R}{\to} \R^2$. The $\Phi$--parallel section at $(0,0,0)$
is $s(x,y)=(x,y,0)$. We have $$\mathfrak{A} \Vdash_{(0,0)} \neg \neg R(x).$$  Let us define $\sigma(t)=(t,t)$, then $\tilde{\sigma}(t)=(t,t,0)$ by equation (\ref{eq: local section and parallel lift}) . It is easy to check that $\sigma^* \mathfrak{A} \nVdash_0 \neg \neg R( \sigma^* s) $. This proves $$\mathfrak{A} \nVdash^{\Phi}_{(0,0,0)} \neg \neg R( x) .$$
\end{exam}
\section{Spatial, horizontal and vertical  extensions}\label{Sec:hor y ver}
As explained in Section~\ref{Sec:epistemological motivation},  an $n$--tuple $(a_1,\cdots,a_n)\in A_m^n$ can be interpreted  as an experimental measurement at a point $m$ of $M$ interpreted as space--time. Intuitively, the sentences  about these measurements should be extensive in space--time and to a certain level independent of the accuracy of the measurements. The horizontal and vertical bundle catch this  distinction between a continuity of the truth depending of the space--time (the spatial or horizontal extension) and a continuity of the truth depending of the accuracy of the measurement (the accuracy extension or vertical extension), respectively.
\begin{defin}
{\upshape (cf. \cite[Secci\'on II]{Caicedo})}
Let $\mathfrak{A}$ be a fiber bundle of $\mathscr{L}$--structures and suppose that $\pi:A \stackrel{F}{\longrightarrow} M$ is endowed with a connection $\Phi$. We define the {\bf spatial extension} of the sections $s_1,\cdots,s_n$ of $A$ for an $\mathscr{L}$--formula $\varphi(x_1, \cdots,x_n)$ in an open subset $U$ of $M$ as follows:
\begin{equation*}
\begin{split}
&[ [\varphi(s_1, \cdots,s_r)] ]_U:=\{u \in U:    \mathfrak{A} \Vdash_u \varphi(s_1, \cdots, s_r) \}.
\end{split}
\end{equation*}
\end{defin}
The next definition generalizes the notion of spatial extension to parallel semantics.
\begin{defin}
Let $\mathfrak{A}$ be a fiber bundle of $\mathscr{L}$--structures and suppose that $\pi:A \stackrel{F}{\longrightarrow} M$ is endowed with a connection $\Phi$. We define the {\bf horizontal extension} of $(e_1,\cdots,e_n )\in A^n_m$ for an $\mathscr{L}$--formula $\varphi(x_1, \cdots,x_n)$ in an open subset $U$ of $M$ as follows:
\begin{equation*}
\begin{split}
&^\Phi[ [\varphi(e_1, \cdots,e_n)] ]_U:=\{u \in U:   \text{ there is a path }\sigma:[0,1]\to U  \text{ such that }\\
&\hspace{0.5cm}\sigma(0)=m, \sigma(1)=u   \text{ and } \mathfrak{A}^\Phi \Vdash_ {\tilde{\sigma}_1(t), \cdots, \tilde{\sigma}_n(t))} \varphi(x_1,\cdots,x_n) \text{ for all } t \in [0,1] \},
\end{split}
\end{equation*}
where $\tilde{\sigma}_i$ are the $\Phi$--lifts of $\sigma$ such that $\tilde{\sigma}_i(0)=e_i$.
\end{defin}
\begin{exam}
Consider the fiber bundle $\pi:\R^2 \stackrel{\R}{\longrightarrow} \R$ (with $\pi(x,y)=x$) with connections $\Phi_1:=dy \otimes\parcial{}{y}$ and $\Phi_2:=(dx+dy)\otimes \parcial{}{y}$. We observe that the vertical bundle of $\pi$ is generated by $\parcial{}{y}$, the $\Phi_1$--horizontal bundle is generated by $\parcial{}{x}$ and the $\Phi_1$--horizontal bundle is generated by  $\parcial{}{x}- \parcial{}{y}$. The identification between the tangent space $T\R$ and the $\Phi_2$--horizontal bundle is given by $d\pi(\parcial{}{x}- \parcial{}{y})=\frac{d}{dt}$. The identification between the tangent space $T\R$ and the $\Phi_1$--horizontal bundle is given by $d\pi(\parcial{}{x})=\frac{d}{dt}$. Given a smooth path $\sigma: (-\epsilon, \epsilon) \to \R$ such that $\sigma(0)=0$, the $\Phi_2$--horizontal lift of $\sigma$ is the solution $\tilde{\sigma}: (-\epsilon, \epsilon) \to \R^2$  of the differential equation 
\begin{center}
$d\pi^{-1}(\sigma'(t)\frac{d}{dt})=\tilde{\sigma}'(t)=(-\sigma'(t),\sigma'(t))$ with initial condition $\tilde{\sigma}(0)=(0,0)$.    
\end{center}
Similarly, the $\Phi_2$--horizontal lift of $\sigma$ is the solution $\tilde{\sigma}: (-\epsilon, \epsilon) \to \R^2$  of the differential equation 
\begin{center}
$d\pi^{-1}(\sigma'(t)\frac{d}{dt})=\tilde{\sigma}'(t)=(\sigma'(t),0)$ with initial condition $\tilde{\sigma}(0)=(0,0)$.    
\end{center} 
If the signature $\mathscr{L}$ has a unary relation symbol $R$, suppose that we have a fibre bundle of $\mathscr{L}$--structures over the fiber bundle $\pi:\R^2 \stackrel{\R}{\longrightarrow} \R$ (with $\pi(x,y)=x$), such that $R^{\mathfrak{A}}:=B_1(0)$, the open  unitary ball  in $\R^2$. 
 With the information given above, we can deduce that the horizontal extensions at $0=(0,0)$ of  $\varphi(x):=R(x)$  associated to the connections $\Phi_1$ and $\Phi_2$ are
$$^{\Phi_1}[[\varphi(0)]]_{\R}=(-1,1) \text{ and } ^{\Phi_2}[[\varphi(0)]]_{\R}=(-\frac{\sqrt{2}}{2},\frac{\sqrt{2}}{2}).$$
\end{exam}
In the previous example we obtained that the horizontal extensions are  open subsets of the base space. This is not always the case.
\begin{exam}
Consider the fiber bundle $\pi:\R^3 \stackrel{\R}{\longrightarrow} \R^2$ (with $\pi(x,y,z)=(x,y)$) with connection $\Phi:=dz \otimes\parcial{}{z}$. We observe the vertical bundle of $\pi$ is generated by $\parcial{}{z}$. The $\Phi$--horizontal bundle is generated by $\parcial{}{x}$ and $\parcial{}{y}$. The identification between the tangent space $T\R^2$ and the $\Phi$--horizontal bundle is given by $d\pi(\parcial{}{x})=\parcial{}{x}$ and $d\pi(\parcial{}{y})=\parcial{}{y}$. Given a smooth path $\sigma: (-\epsilon, \epsilon) \to \R^2$, $\sigma(t):=(\sigma_1(t),\sigma_2(t))$, such that $\sigma(0)=(0,0)$, the $\Phi$--horizontal lift of $\sigma$ is the solution $\tilde{\sigma}: (-\epsilon, \epsilon) \to \R^3$  of the differential equation 
\begin{center}
$d\pi^{-1}(\sigma'(t))=d\pi^{-1}(\sigma'_1(t)\parcial{}{x}+\sigma'_2(t)\parcial{}{y})=\tilde{\sigma}'(t)=(-\sigma'_1(t),\sigma'_2(t),0)$ with initial condition $\tilde{\sigma}(0)=(0,0,0)$.    
\end{center}

If the signature $\mathscr{L}$ has a unary relation symbol $R$, then define $\mathfrak{A}$ a fibre bundle of $\mathscr{L}$--structures such that $R^{\mathfrak{A}}:=\R^3-\{t(1,0,0):t \in \R\}$. We observe that $^\Phi[[\neg R(0,0,0)]]_{\R}=\{(t,0):t \in \R\}$ which is not an open subset. Moreover, $^\Phi[[ R(0,0,0)]]_{\R}$ is the empty set that indicates that there is not an inductive relation in formulas for the horizontal extension of formulas.
\end{exam}
Next we define the vertical extension of a formula. As mentioned before, the intuition of the vertical extension is  how much the validity of the formula  depends of the accuracy of the experimental  measurement. 
\begin{defin}
Let $\mathfrak{A}$ be a fiber bundle of $\mathscr{L}$--structures and suppose that $\pi:A \stackrel{F}{\longrightarrow} M$ is endowed with a connection $\Phi$. We define the {\bf vertical extension} of an $\mathscr{L}$--formula $\varphi(x_1, \cdots,x_n)$ at $(e_1,\cdots,e_n )\in (\oplus_{k=1}^n A)_m$ in an open subset $W$ of the fiber $(\oplus_{k=1}^n A)_m$ as follows:
\begin{equation*}
\begin{split}
^\Phi( (\varphi(e_1, \cdots,e_n)))_W:=&\{(f_1, \cdots,f_n) \in (\oplus_{k=1}^n A)_m: \text{ there is a path }\alpha:[0,1]\to W \\ &\hspace{1cm}\text{ such that } \alpha(0)=(e_1, \cdots,e_n), \alpha(1)=(f_1, \cdots,f_n)\\
&\hspace{1cm} \text{ and }\mathfrak{A}^\Phi \Vdash_ {\alpha(s)} \varphi(x_1,\cdots,x_n) \text{, where } s \in [0,1]\}
    \end{split}
\end{equation*}
\end{defin}
\begin{prop}\label{prop46}
Let  $\mathfrak{A}$ be a fiber bundle of $\mathscr{L}$-structures with underlying fiber bundle $\pi:A\stackrel{F}{\to} M$ whose fiber $F$ and base space $M$ are connected then, for all connection $\Phi$ the vertical extension of the  $\mathscr{L}$-formula $\varphi(x,y):=(x=y)$ at $(e,e) \in (\oplus_{i=1}^2 A)_m$ is $^\Phi(( \varphi ))_M:=\{(f,f):f\in A_m\}$ and the horizontal extension is $[[\varphi]]_M:=M$.
\end{prop}
\bdem
Since $M$ is path--connected, for all $m'\in M$ there is smooth path $\sigma:[0,1] \to M$ such that $\sigma(0)=m$ and $\sigma(1)=m'$ and by uniqueness of the parallel lift we have $\mathfrak{A}^\Phi \Vdash_{(\tilde{\sigma}(s),\tilde{\sigma}(s))}(x=y)$, hence horizontal extension of $(x=y)$ at $(e,e)$ is $M$.
\\
\\
Since $F$  is path--connected,  given $f \in A_m$ there exists a path $\alpha:[0,1]\to A_m$ such that $\alpha(0)=e$ and $\alpha(1)=f$. The path $\beta(s):= (\alpha(s),\alpha(s))$ satisfies $\mathfrak{A}^\Phi \Vdash_{\beta(s)}(x=y)$.
\edem[Prop.~\ref{prop46}]

We believe that the formalization that we offer of horizontal and vertical extensions of a formula could help to  clarify interactions between Geometry, Physics and Mathematical Logic.
\appendix
\section{Geometric background}\label{sec:background}
In this  section we indicate {\it briefly} the basics  of  Differential Geometry needed to understand this article.
\subsection{Fiber bundles and and their connections}
Fiber bundles are spaces that {\it locally  look like}  Cartesian products. They provide a geometric formalization of the idea of  continuous families of  spaces all of which are diffeomorphic.
\begin{defin}\label{defin:fiber bundle}
A $C^k$--fiber bundle ($k\in \N$) consists of three $C^k$--manifolds $A$, $M$ and $F$ and a map $\pi:A \to M$  such that:
\begin{itemize}
\item[i)] $\pi$ is surjective.
\item[ii)]  For each $m \in M$ there exists an open neighborhood $U \subseteq M$ and a $C^k$--diffeomorphism $\psi: \pi^{-1}(U) \to U \times F$ such that the following diagram commutes
$$
\begin{CD}
\pi^{-1}(U) @>\psi>> U \times F\\
@VV{\pi}V @VV{\pi_1}V\\
U @>id>> U
\end{CD}
$$
We denote this fiber bundle by $\pi:A \stackrel{F}{\longrightarrow} M$.
\end{itemize}
We will call $A$ the {\bf total space}, $M$ the {\bf base space} and $F$ the {\bf standard fiber} of the fiber bundle $\pi:A \stackrel{F}{\longrightarrow} M$. The functions $\psi$ in ii) are called the {\bf local trivializations} of the fiber bundle $\pi:A \stackrel{F}{\longrightarrow} M$. A {\bf $C^\infty$--fiber bundle} is a $C^k$--fiber bundle for every $k\in \N$.
\end{defin}
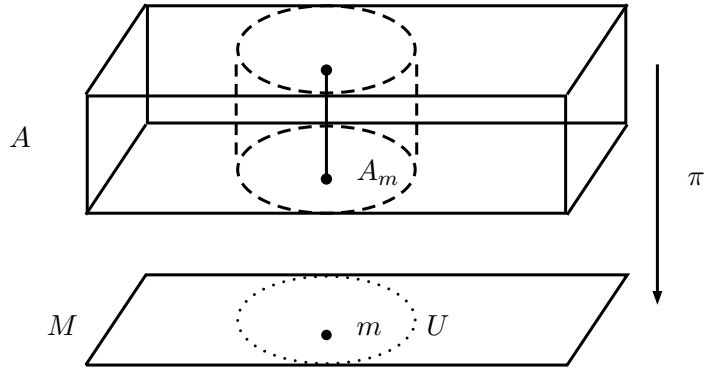
\begin{center}
\begin{pspicture}(0,-2.4155512)(8.226646,2.4155512)
\psframe[linecolor=black, linewidth=0.04, dimen=outer](7.216646,2.4044652)(0.8166464,0.8044653)
\psframe[linecolor=black, linewidth=0.04, dimen=outer](6.4166465,1.2044653)(0.016646422,-0.3955347)
\psline[linecolor=black, linewidth=0.04](0.016646422,1.2044653)(0.8166464,2.4044652)
\psline[linecolor=black, linewidth=0.04](6.4166465,1.2044653)(7.216646,2.4044652)
\psline[linecolor=black, linewidth=0.04](6.4166465,-0.3955347)(7.216646,0.8044653)
\psline[linecolor=black, linewidth=0.04](0.016646422,-0.3955347)(0.8166464,0.8044653)
\psline[linecolor=black, linewidth=0.04](0.016646422,-2.3955348)(0.8166464,-1.1955347)(7.216646,-1.1955347)(6.4166465,-2.3955348)(0.016646422,-2.3955348)
\psline[linecolor=black, linewidth=0.04, arrowsize=0.05291667cm 2.0,arrowlength=1.4,arrowinset=0.0]{->}(7.6166463,1.6044652)(7.6166463,-1.5955347)
\rput[bl](8.016646,0.004465294){$\pi$}
\psellipse[linecolor=black, linewidth=0.04, linestyle=dotted, dotsep=0.10583334cm, dimen=outer](3.2166464,-1.7955346)(1.2,0.6)
\psellipse[linecolor=black, linewidth=0.04, linestyle=dashed, dash=0.17638889cm 0.10583334cm, dimen=outer](3.2166464,0.2044653)(1.2,0.6)
\psellipse[linecolor=black, linewidth=0.04, linestyle=dashed, dash=0.17638889cm 0.10583334cm, dimen=outer](3.2166464,1.8044653)(1.2,0.6)
\psline[linecolor=black, linewidth=0.04, linestyle=dashed, dash=0.17638889cm 0.10583334cm](2.0166464,0.4044653)(2.0166464,1.6044652)
\psline[linecolor=black, linewidth=0.04, linestyle=dashed, dash=0.17638889cm 0.10583334cm](4.4166465,0.4044653)(4.4166465,1.6044652)
\psline[linecolor=black, linewidth=0.04, dotsize=0.07055555cm 2.0]{**-**}(3.2166464,0.004465294)(3.2166464,1.6044652)
\psdots[linecolor=black, dotsize=0.14](3.2166464,-1.9955347)
\rput[bl](3.6166465,0.004465294){$A_m$}
\rput[bl](-1,0.5){$A$}
\rput[bl](-0.5,-2){$M$}
\rput[bl](3.6166465,-1.9955347){$m$}
\rput[bl](4.4166465,-1.9955347){
$U$}
\end{pspicture}
\end{center}
Fiber bundles have different notions of morphism which provide different categories. We will be interested in morphisms between fiber bundles over the same base.
\begin{defin}\label{defin: morphism fiber bundle}
A {\bf morphism of fiber bundles} between $\pi_1:A_1 \to M$ and $\pi_1:A_1 \to M$ is a  smooth map  $\Phi:A_1 \to A_2$ such that the following diagram commutes
$$
\begin{CD}
A_1 @>\Phi>>A_2\\
@VV{\pi_1}V @VV{\pi_2}V\\
M @>i>> M
\end{CD}
$$
where the map $i:M \to M$ is the identity.
\end{defin}
Throughout this article, we assume that all the fiber bundles are $C^\infty$.
\begin{defin}
\begin{itemize}
\item A global section of a fiber bundle $\pi:A \stackrel{F}{\longrightarrow} M$ is a $C^\infty$--function $s:M\to A$ such that $\pi\circ s(m)=m $ for all $m\in M$. The set of all global sections of $\pi:A \stackrel{F}{\longrightarrow} M$ will be denoted by $\Gamma(A)$.
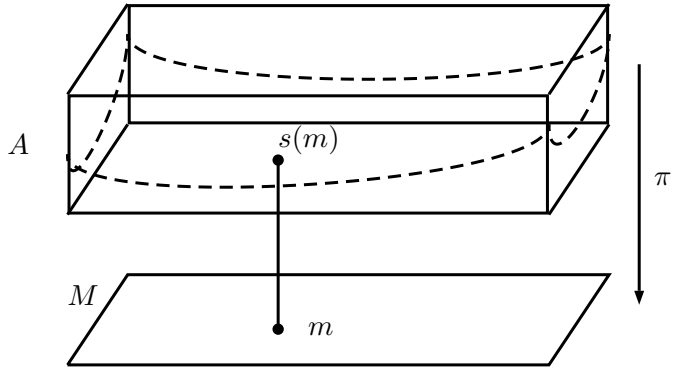
\begin{center}
\psscalebox{1.0 1.0} 
{
\begin{pspicture}(0,-2.6100082)(9.01,2.6100082)
\psframe[linecolor=black, linewidth=0.04, dimen=outer](8.4,2.2100081)(2.0,0.61000824)
\psframe[linecolor=black, linewidth=0.04, dimen=outer](7.6,1.0100082)(1.2,-0.58999175)
\psline[linecolor=black, linewidth=0.04](1.2,1.0100082)(2.0,2.2100081)
\psline[linecolor=black, linewidth=0.04](7.6,1.0100082)(8.4,2.2100081)
\psline[linecolor=black, linewidth=0.04](7.6,-0.58999175)(8.4,0.61000824)
\psline[linecolor=black, linewidth=0.04](1.2,-0.58999175)(2.0,0.61000824)
\psline[linecolor=black, linewidth=0.04](1.2,-2.5899918)(2.0,-1.3899918)(8.4,-1.3899918)(7.6,-2.5899918)(1.2,-2.5899918)
\psline[linecolor=black, linewidth=0.04, arrowsize=0.05291667cm 2.0,arrowlength=1.4,arrowinset=0.0]{->}(8.8,1.4100082)(8.8,-1.7899917)
\pscustom[linecolor=black, linewidth=0.04]
{
\newpath
\moveto(5.6,-2.1899917)
}
\pscustom[linecolor=black, linewidth=0.04]
{
\newpath
\moveto(9.6,-2.5899918)
}
\pscustom[linecolor=black, linewidth=0.04]
{
\newpath
\moveto(1.2,-4.5899916)
}
\pscustom[linecolor=black, linewidth=0.04]
{
\newpath
\moveto(3.2,-4.189992)
}
\pscustom[linecolor=black, linewidth=0.04]
{
\newpath
\moveto(6.8,-1.7899917)
}
\psbezier[linecolor=black, linewidth=0.04, linestyle=dashed, dash=0.17638889cm 0.10583334cm](1.2,0.21000823)(1.2,-0.58999175)(2.0,1.0100082)(2.0,1.8100082397460937)
\psbezier[linecolor=black, linewidth=0.04, linestyle=dashed, dash=0.17638889cm 0.10583334cm](7.6,0.61000824)(7.6,-0.18999176)(8.4,1.0100082)(8.4,1.8100082397460937)
\psbezier[linecolor=black, linewidth=0.04, linestyle=dashed, dash=0.17638889cm 0.10583334cm](2.0,1.8100083)(2.0,1.0100082)(8.4,1.0100082)(8.4,1.8100082397460937)
\psbezier[linecolor=black, linewidth=0.04, linestyle=dashed, dash=0.17638889cm 0.10583334cm](1.2,0.21000823)(1.2,-0.58999175)(7.6,-0.18999176)(7.6,0.6100082397460938)
\rput[bl](9,-0.18999176){$\pi$}
\psline[linecolor=black, linewidth=0.04, dotsize=0.07055555cm 2.0]{**-**}(4.0,-2.1899917)(4.0,0.21000823)
\rput[bl](4.4,-2.1899917){$m$}
\rput[bl](1.2,-1.7899917){$M$}
\rput[bl](0.4,0.21000823){$A$}
\rput[bl](4.0,0.21000823){$s(m)$}
\end{pspicture}
}
\end{center}

\item A local section of a fiber bundle $\pi:A \stackrel{F}{\longrightarrow} M$ is a $C^\infty$--function $s:U\to A$ defined on a open neighborhood $U$ of $M$ such that $\pi\circ s(m)=m $ for all $m\in U$
\end{itemize}
\end{defin}
\begin{rem}\label{Rem:TangetSpace}
If $M$ is a manifold, it is a usual exercise in differential geometry to observe that $TM:=\bigcup_{m \in M} T_mM$ (the disjoint union of all tangent spaces of $M$)  is the total space of a fiber bundle with basis $M$ and fiber $\R^n$  where $n$ is the dimension of the manifold $M$ and $\pi:TM \to M$ is the natural projection. Special types of fiber bundles are the the principal bundles and vector bundles.
\end{rem}
\begin{notat}
Let  $f:M \to N$ be a smooth function between two manifolds $M$ and $N$. We will denote the derivative of $f$ in $m\in M$ by $df_m:T_mM \to T_{f(m)}N$.
\end{notat}
Let $\pi:A \stackrel{F}{\longrightarrow} M$ be a  fiber bundle. For each $a \in A$, the tangent space $T_a (A_m)$ of the fiber at $m:=\pi(a)$ define a vector space of directions in $T_aA$ called {\bf vertical directions}, more formally:
\begin{defin}
The subvector bundle  $VA:={ \rm Ker} (d\pi)$ of $TA$ is called {\bf the vertical bundle} of the fiber bundle $\pi:A \stackrel{F}{\longrightarrow} M$.
\end{defin}
We observe that $T_a(A_m)=(VA)_a$ where $a\in A_m$. To define horizontal directions we need a  connection   on $\pi:A \stackrel{F}{\longrightarrow} M$, that is basically a choice of a projection $\Phi_a: T_aA \to T_aA$ on $(VA)_a$ for each $a\in A_m$.
\begin{defin}\label{Connection}
A {\bf connection} $\Phi$ for a fiber bundle $\pi:A \stackrel{F}{\longrightarrow} M$ is a smooth $1$--form  of $A$ with values  in $VA$ such that  for each $a \in A$,  $\Phi_a^2=\Phi_a$ and ${\rm Im }(\Phi_a)=(VA)_a$. \end{defin}
A connection $\Phi$ belongs to $\Omega^1(A,VA)$ and   $\Phi(a)$ can be thought as a linear map from $T_aA $ to $T_a A$ for each $a \in A$.
A connection $\Phi$ induces a notion of {\bf horizontal bundle $HA_{\Phi}$}
\begin{equation}\label{eq:definition HA}
HA_{\Phi}:=HA:={\rm Im}(I_{TA}-\Phi) \subseteq TA.    
\end{equation}
Moreover, for all $a \in A$, we have that $d\pi_a:(HA_{\Phi})_a \to T_{\pi(a)} M $ is a  canonical isomorphism, and it is easy to see that $TA=VA\oplus HA$. In particular, the tangent directions of the total space $A$ are decomposed in horizontal and vertical directions.
\\
\\
The parallel transport of a curve $\sigma:(-r,r) \to M$ at a point $a \in A$ is the lift $\tilde{\sigma}: (-r,r) \to A$ of $\sigma$ (i.e $\pi \circ \tilde{\sigma}=\sigma$) whose velocity belongs to the horizontal direction and such that $\tilde{\sigma}(0)=a$.   The following theorem formalizes this notion and guarantees that, for every connection $\Phi$ of $A$ and every curve $\sigma$ and point $a \in A$, there exists (locally) a unique parallel transport.
\begin{thm}\label{thm: existence parallel transport}{\upshape \cite[Theorem 9.8]{KolarMichor}} Let $\pi:A \stackrel{F}{\longrightarrow} M$ be a fiber bundle with connection $\Phi$ and let $\sigma:(-r,r ) \to M$ be a smooth curve such that $\sigma(0)=m$. Then, there exists a neighborhood $U$ of $A_m \times \{0\}$ in  $A_m \times (-r,r)$ and a smooth function $\tilde{\sigma}:U \to A$ such that:
\begin{itemize}
\item[i)] $\pi(\tilde{\sigma}(a,t))=\sigma(t)$ for all $(a,t) \in U \subseteq A_m \times (-r,r)$ and $\tilde{\sigma}(a,0)=a$.
\item[ii)] $\Phi(\parcial{}{t}\tilde{\sigma}(a,t) )=0$ for all $(a,t) \in U$.
\item[iii)] $U$ is maximal with respect to i) and ii).
\end{itemize}
\end{thm}
With the notation of the previous theorem, let us recall that $d\pi_{a}$ is an isomorphism between the horizontal bundle $HA_{a}$ at $a$ of the connection $\Phi$ (see (\ref{eq:definition HA})) and the tangent space $T_aA$. With this identification, condition ii) of Theorem~\ref{thm: existence parallel transport} can be written 
\begin{equation}\label{eq: parallel transform horizontal version}
d\pi^{-1}_{\sigma(t)}(\sigma'(t))=\tilde{\sigma}'(t).
\end{equation}
Using the notation of Theorem~\ref{thm: existence parallel transport}.
\begin{defin}\label{Def:lift of a path}
\begin{itemize}
\item[i)] Given $a\in A_m$ ($m\in M$),  the function $t \mapsto\tilde{\sigma}(a, t)$ defined in the previous theorem is called {\bf parallel transport} along the curve $\sigma$ of $a$ (associated to the connection $\Phi$).
\item[ii)] A connection $\Phi$ on $\pi:A \stackrel{F}{\longrightarrow} M$ is called a {\bf complete connection}, if the parallel transport $\tilde{\sigma}$ along any smooth curve $\sigma : (-r, r) \to M$ is defined in all elements belonging to $A_{\sigma(0)}\times (-r, r)$.
\end{itemize}
Also we call $t \mapsto\tilde{\sigma}(a, t)$  the {\bf horizontal lift} of $\sigma$ at $a$.
\end{defin}
 Intuitively, parallel transport formalizes the notion of a movement on a configuration space that does not change the internal states (see Section~\ref{Sec:epistemological motivation}). The notion of {\it completeness of a connection} is a technical condition which will simplify this presentation. Complete connections are also called {\it Ehresmann connections}. 
The following theorem allows us to consider a complete connection in any fiber bundle, which helps us to avoid technicalities. This is the reason because we assume completeness of all connections considered along this article.

\begin{thm}{\upshape \cite[Page 81]{KolarMichor}}
Each fiber bundle admits complete connections.
\end{thm}
Next, we define the notion of pullback of a fibre bundle.
\begin{defin}\label{Def:Pullback fiber bundle}
Given a smooth function $f:N \to M$ the {\bf pullback} of a fiber bundle $A$ is the fiber bundle $f^*(A)$ whose total space is  $$f^*(A):=\{(n,a)\in N \times A: f(n)=\pi(a)\}$$ with  projection $\pi_{f^*A}(n,a)=n$, the topology inherited from $N \times A$ and differential structure naturally defined from the trivializations induced by the fiber bundle $A$.
\end{defin}
Given  a smooth function $f:N \to M$ there is a natural morphism of fiber bundles $\tilde{f}:f^*(A) \to A$   defined by $\tilde{f}(n,a)=a$ and illustrated by the following diagram. 
\begin{center}
\begin{tikzpicture}[x=0.75pt,y=0.75pt,yscale=-1,xscale=1]
\draw    (110.5,137) -- (169.5,137) ;
\draw [shift={(171.5,137)}, rotate = 180] [color={rgb, 255:red, 0; green, 0; blue, 0 }  ][line width=0.75]    (10.93,-3.29) .. controls (6.95,-1.4) and (3.31,-0.3) .. (0,0) .. controls (3.31,0.3) and (6.95,1.4) .. (10.93,3.29)   ;
\draw    (196.5,76) -- (196.5,124) ;
\draw [shift={(196.5,126)}, rotate = 270] [color={rgb, 255:red, 0; green, 0; blue, 0 }  ][line width=0.75]    (10.93,-3.29) .. controls (6.95,-1.4) and (3.31,-0.3) .. (0,0) .. controls (3.31,0.3) and (6.95,1.4) .. (10.93,3.29)   ;
\draw  [dash pattern={on 4.5pt off 4.5pt}]  (125,63) -- (179.5,63.96) ;
\draw [shift={(181.5,64)}, rotate = 181.01] [color={rgb, 255:red, 0; green, 0; blue, 0 }  ][line width=0.75]    (10.93,-3.29) .. controls (6.95,-1.4) and (3.31,-0.3) .. (0,0) .. controls (3.31,0.3) and (6.95,1.4) .. (10.93,3.29)   ;
\draw  [dash pattern={on 4.5pt off 4.5pt}]  (96.5,80) -- (96.5,118) ;
\draw [shift={(96.5,120)}, rotate = 270] [color={rgb, 255:red, 0; green, 0; blue, 0 }  ][line width=0.75]    (10.93,-3.29) .. controls (6.95,-1.4) and (3.31,-0.3) .. (0,0) .. controls (3.31,0.3) and (6.95,1.4) .. (10.93,3.29)   ;
\draw (99,59) node  [align=left] {$\displaystyle f^{*}( A)$};
\draw (198,62) node  [align=left] {$\displaystyle A$};
\draw (94,137) node  [align=left] {$\displaystyle N$};
\draw (195,137) node  [align=left] {$\displaystyle M$};
\draw (215,102) node  [align=left] {$\displaystyle \pi $};
\draw (149,42) node  [align=left] {$\displaystyle \tilde{f}$};
\draw (77,101) node  [align=left] {$\displaystyle \pi_{N}$};
\draw (139,154) node  [align=left] {$\displaystyle f$};
\end{tikzpicture}
\end{center}
Given a fiber bundle $\pi: A \stackrel{F}{\longrightarrow} M$ with fiber $F$ and $k\in \N\setminus \{0\}$, we  define its {\bf $k$--direct sum} (denoted by denote by $\oplus_{i=1}^ k A$) as  a fiber bundle over $M$ whose fiber at $m \in M$ is given by $A^k_m$.  In contrast, the {\bf Cartesian power} $A^k$ corresponds to the fiber bundle over the Cartesian power $M^k$ with the natural projection $\pi \times \cdots \times \pi$.
\begin{defin}\label{DirectSum}
 Let   $d:M \to \Pi_{i=1}^k M$  be the diagonal function defined by $d(m):=(m,m, \cdots, m)$. The {\bf $k$--direct sum of} a fiber bundle $\pi:A \stackrel{F}{\longrightarrow} M$,
denoted by $\oplus_{i=1}^ k A$, is defined by
$$
\oplus_{i=1}^ k A:=d^*(A^k).
$$
\end{defin}
From a categorical point of view the direct sum defined in Definition~\ref{DirectSum} is the product of the category of fiber bundles over a fixed basis. In the categories of vector spaces and vector bundles the direct sum coincide with the product, we hope this justify the abuse of terminology. 
\\
\\
We can also pullback sections of fiber bundles to sections of (the pullbacked) fiber fundle.
\begin{defin}\label{Def:Pullback connection}
Given a fiber bundle $\pi:A \to M$ and a function $f:N \to M$, we define {\bf the pullback of a section $s:U \to M$} as the section $f^*s$ of $f^*A$ defined by
$$
f^*s(n):=(n,s(f(n))).
$$
\end{defin}

The next proposition claims that it is possible to pullback connections.
\begin{prop}\label{Prop:pullback of connections}
Let $\Phi$ be a connection on  the  fiber bundle $\pi:A \to M$ and let $f:N \to M$ be a smooth function. $f^*(\Phi)$ induces a  connection on $f^*(A)$. 
\end{prop}
\bdem
$f^*(\Phi)$ is a 1--form of $f^*(A)$ with values in $VA$, hence we can think of it as an homomorphism from $Tf^*(A)$ to $VA$. Recall the diagram below Definition~\ref{Def:Pullback fiber bundle} for the definition of $\tilde{f}$.  $d\tilde{f}$ induces an identification of the vertical bundle of the pullback $Vf^*(A)$ and the vertical bundle of $A$. Hence we can define a connection $\varphi$  on $f^*(A)$ by $\varphi_{(n,a)}:=df^{-1}_n\circ f^*(\Phi)_{(n,a)}$. Using that  $f^*(\Phi)_{(n,a)}:T_{(n,a)}f^*(A) \to VA$ is given by $\Phi_a \circ df_n $, we can observe that $\varphi_{(n,a)}:T_{(n,a)}f^*(A) \to T_{(n,a)}f^*(A)$ is a projection on $V_{(n,a)}f^*(A)$.\edem[Proposition~\ref{Prop:pullback of connections}]
With an slight abuse of  notation, we will denote the connection $\varphi$ in the previous proposition by $f^*(\Phi)$.
\begin{cor}\label{Cor:pullback connection dir sum}
A connection $\Phi$ on $\pi:A \stackrel{F}{\longrightarrow} M$ induces a connection on  the $k$--direct sum $\oplus_{i=1}^ k A $  that we will denote by $\oplus_{i=1}^ k \Phi$.
\end{cor} 
\bdem
The connection $\Phi$ on the fiber bundle $\pi:A \to M$ induces the connection $\Pi_{i=1}^k \Phi$    on the 
Cartesian product $\Pi_{i=1}^k A$ (thought as a fiber bundle over $\Pi_{i=1}^k M$). This connection is defined by $\Pi_{i=1}^k \Phi_{(e_1,\cdots,e_k)}=\Phi_{e_1}\oplus \cdots \oplus \Phi_{e_k}$, where we are identifying the tangent space $T\Pi_{i=1}^k A_{(e_1,\cdots,e_k)}$ naturally with  $\oplus_{i=1}^k T_{e_i}A$. This identification induces a connection on $\oplus_{i=1}^k A$ taking $d^*(\Pi_{i=1}^k \Phi)$ (denoted by $\oplus_{i=1}^k \Phi$) where $d:M \to \Pi_{i=1}^k M$ is the diagonal map.
\edem[Corollary \ref{Cor:pullback connection dir sum}]
\begin{prop}\label{Prop:lift of Ak}
Let $\Phi$ be a connection on  the  fiber bundle $\pi:A \to M$ and let $\sigma:(-\varepsilon,\varepsilon) \to M$ be a path such that $\sigma (0)=m$. Let $(m,e_1,\cdots,e_k) \in \oplus_{i=1}^k A$ and let us denote $\tilde{\sigma}_i$ the $\phi$--lift of $\sigma$ to the fiber bundle $A$ such that $\tilde{\sigma}_i(0)=e_i$. If $\alpha$ is the $\oplus_{i=1}^k \Phi$--lift of $\sigma$ to the fiber bundle $\oplus_{i=1}^k A$ such that $\alpha(0)=(m,e_1,\cdots,e_k)$, then $\alpha(t)=(\sigma(t),\tilde{\sigma}_1(t),\cdots,\tilde{\sigma}_k(t))$. 
\end{prop}
\bdem
We observe that the inclusion of $\oplus_{i=1}^k A$ in  $A \times \cdots \times A$ is an embedding. We can see that $H_{\oplus_{i=1}^k \Phi}$ the horizontal bundle of  $\oplus_{i=1}^k A$  is equal to 
$$(H_{\oplus_{i=1}^k \Phi})_{(a_1,\cdots,a_k)}=$$$$\{ (v_1,\cdots,v_k)\in T(A \times \cdots \times A)_{a_1,\cdots a_k}):d\pi_{a_1}(v_1)=\cdots=d\pi_{a_k}(v_k)\},$$
where $\pi:A \to M$ denotes the projection of the fibre bundle $A$. With this identification, for all $t \in (-\varepsilon,\varepsilon)$, $\alpha'(t)=(\tilde{\sigma}_1'(t),\cdots,\tilde{\sigma}'_k(t))$ since $\tilde{\sigma}_i'(t)=d\pi^{-1}(\sigma'(t))$ for $d\pi^{-1}$ restricted from $TM$ to the $\Phi$--horizontal bundle $HA$. Since $\alpha(0)=(\tilde{\sigma}_1(0),\cdots,\tilde{\sigma}_k(0))$, we have proved the proposition. \edem[Proposition~\ref{Prop:lift of Ak}]
\ \\
\ \\
\begin{prop}\label{prop: horizontal lift pullback}
Let $f:M \to N$  be a smooth function and let $\pi:A \to N$ be a fiber bundle with fiber $F$ and connection $\Phi$. Suppose that $\sigma:(-\epsilon,\epsilon ) \to M$ is a smooth path such that $\sigma(0)=m$ and let $(m,e) \in f^*(A)_m \subset M \times A$ be fixed. Let us denote by $\tilde{\alpha}$ the $\Phi$--horizontal lift of $\alpha:=f\circ \sigma$
such that $\tilde{\alpha}(0)=e$. Then $\tilde{\sigma}$, the $f^* \Phi$--horizontal lift of $\sigma$ such that  $\tilde{\sigma}(0)=(m,e)$, is equal to $ t \mapsto (\sigma(t),\tilde{\alpha}(t))$.
\end{prop}
\bdem
Let us denote $\beta(t):=(\sigma(t),\tilde{\alpha}(t))$, clearly $\pi_{f^* A}\beta=\sigma(t)$. We have $\beta '(t)=(\sigma '(t),\tilde{\alpha}'(t))$. We recall that $f^*\Phi=d\tilde{f}^{-1} \Phi \tilde{f}$ (see Proposition~\ref{Prop:pullback of connections}). Then $f^*(\Phi)\beta '(t)=d\tilde{f}^{-1} \Phi \tilde{f} (\sigma '(t),\tilde{\alpha}'(t))= d\tilde{f}^{-1} \Phi \tilde{\alpha}'(t) =0$. 
\edem[Propostion~\ref{prop: horizontal lift pullback}]
Next proposition synthetizes some canonical equivalences of the pullback which are used along the article.
\begin{prop}\label{Prop:canonical isomorphisms}
Let $f:M \to N$ and $g:N \to P$ be smooth functions and let $\pi:A \to P$ be a fiber bundle with fiber $F$. Then, 
\begin{itemize}
    \item[i)] The fiber bundle $(g \circ f)^*(A)$ is canonically isomorphic to $f^*(g^*(A))$.
    \item[ii)] The fiber bundle $\oplus_{i=1}^k f^*(A)$ is canonically isomorphic $f^*( \oplus_{i=1}^k A)$.  
    \item[iii)] Let $U$ be an open subset of $P$, and let us denote $\Gamma(U,A)$ the sections of the fiber bundle $A$ whose domain contains  $U$. We have that there is a canonical identification of $\Gamma(U,\oplus_{i=1}^k A)$ and  $\Gamma(U, A)^k$.   
\end{itemize}
\end{prop}
\bdem
Let us observe that $(g \circ f)^*(A)=\{(m,a) \in M \times A: (g \circ f)(m)= \pi(a) \}$ and
\begin{equation*}
\begin{split}
f^* (g^*(A))&=\{(m,b)\in M \times g^*(A):f(m)=\pi_{f^*(A)}(b) \}
\\
&=\{(m,n,a)\in M \times N \times A: f(m)=n \text{ and } g(n)=\pi(a) \}.      
\end{split}
\end{equation*}
It is straightforward to see that the canonical isomorphism is given by $(m,a) \mapsto (m,f(m),a)$.
\\
\\
We have that 
$$
f^*( \oplus_{i=1}^k A)=\{(n,m,a_1,\cdots,a_k)):f(n)=m \text{ and } a_i \in A_m\}
$$
and
$$\oplus_{i=1}^k f^*(A)=\{(n,((n,a_1),\cdots (n,a_k)) \in N \times f^*(A)^k: a_i \in A \text{ and } \pi(a_i)=f(n)\}.$$
The isomorphism between $f^*( \oplus_{i=1}^k A)$ and $\oplus_{i=1}^k f^*(A)$ is given by $(n,((n,a_1),\cdots (n,a_k)) \mapsto (n,f(n),a_1,\cdots,a_k)$.
\\
\\
The identification of $\Gamma(U,\oplus_{i=1}^k A)$ and  $\Gamma(U, A)^k$ is given by $(s_1,\cdots, s_k) \mapsto \vec{s}$ where $\vec{s}(p)=(p,s_1(p),\cdots, s_k(p))$ for $p \in P$. Given  $\vec{s} \in \Gamma(U,\oplus_{i=1}^k A)$, with $\vec{s}(p)=(p,s_1(p),\cdots, s_k(p))$,  its inverse under this map is the $k$--tuple $(s_1,\cdots, s_k) \in \Gamma(U, A)^k$.
\edem[Proposition~\ref{Prop:canonical isomorphisms}]
\ \\
\ \\
Since we work with pullbacks of fiber bundles, sections, and connections through paths, we remark some facts about this.

\begin{prop}\label{Prop: lift of identity as  a section}
Let $\pi:A \to N$ be a fiber bundle with fiber $F$ with a connection $\Phi$. Let $\sigma:(-\epsilon, \epsilon) \to M$ be a path such that $\sigma(0)=m$ and let us denote $\alpha(t):=t$, the identity path on $\R$. Then, we can naturally interpret the $\sigma^*\Phi$--horizontal lifts of $\alpha$, as a section of the fiber bundle $\sigma^*A$. Moreover, given  $e \in A_m$, if $\tilde{\sigma}$ denotes the $\Phi$--horizontal lift of $\sigma$ such that $\tilde{\sigma}(0)=e$  and $s$ is the $\sigma^*\Phi$--horizontal lift of $\alpha$ such that $s(0)=(0,e)$, then $s(t)=(t,\tilde{\sigma}(t))$
\end{prop}
\bdem
$s(t)=(t,\tilde{\sigma}(t))$ is a direct consequence of Proposition~\ref{prop: horizontal lift pullback}. Clearly, $s(t)$ is a section of $\sigma^* A$. \edem[Proposition~\ref{Prop: lift of identity as  a section}]

\begin{defin}\label{defin: isomorphic bundles with connection}
Let $\pi_1:A_1 \to M$ and $\pi_2:A_2 \to M$ be fiber bundles with connections $\Phi_1$ and $\Phi_2$. We say that these {\bf fiber bundles are isomorphic} if there is an isomorphism of fiber bundles
$$
\begin{CD}
A_1 @>h>>A_2\\
@VV{\pi_1}V @VV{\pi_2}V\\
M @>i>> M
\end{CD}
$$
which preserves the horizontal bundles induced by the connection,  explicitly $dh (HA_1)=HA_2$. 
\end{defin}

Isomorphism of fiber bundles with connections send horizontal liftings in horizontal liftings.

As a consequence of the previous propositions, we have.
\begin{cor}
Let $f:M \to N$ and $g:N \to P$  be  smooth functions and let $\pi:A \to N$ be a fiber bundle with fiber $F$ and connection $\Phi$. The following pairs of fiber bundles with connection  are isomorphic,
\begin{itemize}
    \item[i)] $(f\circ g)^*A, (f\circ g)^*\Phi)$ and $f^*(g^*A), f^*(g^*\Phi))$.
    \item[ii)] $(f^*(\oplus_{i=1}^k A), f^*(\oplus_{i=1}^ k \Phi))$ and $(\oplus_{i=1}^k f^*A, \oplus_{i=1}^ k f^*\Phi)$. 
\end{itemize}
\end{cor}
We will denote by $$\chi:=I_{TA}-\Phi$$ the projection on the horizontal bundle of a connection $\Phi$.
\begin{defin}\label{curvature}
Let $\Phi$ be a connection on the fiber bundle $\pi:A \stackrel{F}{\longrightarrow} M$.  The {\bf curvature} of $\Phi$ is the two form $R$ with values in $VA$ defined by
$$
R(X,Y):=\Phi([\chi(X),\chi(Y)]),
$$
where $X,Y$ are vector fields of $A$.
\end{defin}
The curvature of the connection is an obstruction (via the Frobenius theorem, see \cite{Warner}) to the integrability of the differential distribution $HA_{\Phi}$ on $A$.
\begin{defin}\label{Def:parallel section}
Given $\Phi$ a connection on $\pi:A\stackrel{F}{\longrightarrow} M$, we will say that {\bf a local section $s:M \to A$ is $\Phi$--parallel} if  the horizontal bundle  $HA_\Phi$ restricted to $U$ is equal to the  image of $ds$ restricted to $TU$.
\end{defin}
We observe that  $\Phi$--parallel sections exist only when the connection $\Phi$ has curvature $0$.
\begin{prop}\label{prop: curvature 0 sections}
Let $\Phi$ be a connection on the fiber bundle $\pi:A \stackrel{F}{\longrightarrow} M$ whose curvature $R$ is $0$. Then, for every $m \in M$, and every $e \in A_m$, there is a neighborhood $U$ of $m$ and a  local section $s:U \to A$ such that $s$ is a $\Phi$--parallel section and $s(m)=e$.
\end{prop}
\bdem
Let us take $m \in M$ and  $e \in A_m$. Since the curvature of $\Phi$ is zero, the horizontal bundle $HA_{\Phi}$ of $\Phi$ is integrable.  The theorem of Frobenius (see \cite{Warner}) implies that there exists an integral submanifold $N$ of $M$ which contains $a$ i.e. a submanifold $N$ such that for all  $n\in N$, $T_nN=(HA_{\Phi})_n$. If we take $f:=\pi \vert_{N}$ then $d f_n:(HA_{\Phi})_a \to T_{\pi(a)}M$ is an isomorphism. The inverse function theorem implies that there exists an open neighborhood $U$ of $m$ such that $f^{-1}:U \to N$ exists. We can take $s:=f^{-1}$.  
\edem[Proposition~\ref{prop: curvature 0 sections}]
\begin{remark}
As a general principle, we work in this article around small open neighborhoods of the point $m \in M$. Hence, the section whose existence is claimed by Proposition~\ref{prop: curvature 0 sections} is unique in the sense that,  if $s_1$ and $s_2$ with connected domains $U_1 \subset M$ and $U_2 \subset M$ respectively, are parallel sections such that $s_1(m)=s_2(m)=e$, then $s_1=s_2$ in the open neighborhood of $m$ $U_1\cap U_2$.  
\end{remark}
Each section $s$ of the previous proposition is a local section of $A$ and  for all path $\sigma:(-\varepsilon,\varepsilon) \to M$ such that $\sigma(0)=m$, we have
\begin{equation}\label{eq: local section and parallel lift}
\tilde{\sigma}(t)=s\circ \sigma(t)    
\end{equation}
 where $\tilde{\sigma}$ is the $\Phi$--lift of a path $\sigma$ such that $\tilde{\sigma}(0)=e$ (see (\ref{eq: parallel transform horizontal version})). 
\begin{rem}\label{rem: identificacion tangente Rn y Rn}
On $\R^n$ there is a natural identification of the tangent space of $\R^n$  $T\R^n$ and $\R^n\times \R^n$, specifically we identify $T_{(p_1,\cdots,p_n)}\R^n$ with $\R^n$ using the linear map $(a_1, \cdots,a_n) \mapsto \sum_{i=1}^n a_i \partial_{x_i}\vert_{(p_1,\cdots,p_n)}$.
\end{rem}
Many examples of this article will be related with the following basic fiber bundle with connection.  
\begin{exam}\label{Exam:trivial fiber bundle with connect}
Consider the fiber bundle $\rho:\R^3 \to \R^2$ with $\rho(x,y,z):=(x,y)$ whose fiber is $\R$. Let us define $\Phi= dz \otimes \partial_z$. Using the identification of Remark~\ref{rem: identificacion tangente Rn y Rn}, we can write $\Phi_{(x,y,z)^t}=\begin{pmatrix}0&0&0\\0&0&0\\0&0&1
\end{pmatrix}$. One can check  that $HA_{(x,y,z)^t}={\rm Span}\{\partial_x \vert_{(x,y,z)^t}, \partial_y \vert_{(x,y,z)^t}\}$. With the identification of Remark~\ref{rem: identificacion tangente Rn y Rn}, $HA_{(x,y,z)^t}= {\rm Span}\left\{\begin{pmatrix}1\\0\\0
\end{pmatrix},\begin{pmatrix}0\\0\\1
\end{pmatrix}\right\}$. The isomorphism of (\ref{eq: parallel transform horizontal version}) is given by $d\pi^{-1}(a\partial_x + b\partial_y)=a\partial_x + b\partial_y$. Let us denote $\sigma(t):=(mt,nt)$, for fixed $m,n$ in $\N$, with this considerations, it is straightforward to see that $\tilde{\sigma}$  the $\Phi$--horizontal lift of $\sigma$ such that $\tilde{\sigma}(0)=(0,0,a)$ is $\tilde{\sigma}(t)=(mt,nt,a)$.
\end{exam}
\section*{Aknowledgments}
The first author is grateful with Andr\'es Villaveces for many inspiring conversations around these topics. He is also grateful with  Xavier Caicedo and Fernando Zalamea whose comments improve the presentation of the article. The second author thanks the first author for the invitation to participate in this project.
\bibliographystyle{alpha}

\bibliography{literatur}

\end{document}